\documentclass[11pt]{article}
\usepackage{amscd}
\usepackage{amsmath}
\usepackage{latexsym}
\usepackage{amsfonts}
\usepackage{amssymb}
\usepackage{amsthm}
\usepackage{graphicx}
\usepackage{verbatim}
\usepackage{mathrsfs}
\usepackage{enumerate}
\usepackage{hyperref}

\theoremstyle{plain}\newtheorem{definition}{Definition}[section]
\theoremstyle{definition}\newtheorem{theorem}{Theorem}[section]
\theoremstyle{plain}\newtheorem{lemma}[theorem]{Lemma}
\theoremstyle{plain}
\theoremstyle{plain}\newtheorem{proposition}[theorem]{Proposition}
\theoremstyle{remark}\newtheorem{remark}{Remark}[section]

\newcommand{\norm}[1]{\left\|#1\right\|}

\begin{document}
\title{On the regularity of a class of generalized quasi-geostrophic equations}
\author{Changxing Miao\footnote{Institute of Applied Physics and Computational Mathematics, P.O. Box 8009,
            Beijing 100088, P.R. China. Email: miao\_{}changxing@iapcm.ac.cn.}\;  and  Liutang Xue\footnote{
            The Graduate School of China Academy of Engineering Physics, P.O. Box 2101, Beijing 100088, P.R. China. Email: xue\_{}lt@163.com.}}

\date{}
\maketitle

\begin{abstract}
  In this article we consider the following generalized quasi-geostrophic equation
  \begin{equation*}
    \partial_t\theta + u\cdot\nabla \theta + \nu \Lambda^\beta \theta =0, \quad u= \Lambda^\alpha \mathcal{R}^\bot\theta, \quad  x\in\mathbb{R}^2,
  \end{equation*}
  where $\nu>0$, $\Lambda:=\sqrt{-\Delta}$, $\alpha\in ]0,1[$ and $\beta\in ]0,2[$. We first show a general criterion yielding the nonlocal maximum principles for
  the whole space active scalars, then mainly by applying the general criterion, for the case $\alpha\in]0,1[$ and $\beta\in ]\alpha+1,2]$ we obtain the global
  well-posedness of the system with smooth initial data; and for the case $\alpha\in ]0,1[$ and $\beta\in ]2\alpha,\alpha+1]$
  we prove the local smoothness and the eventual regularity of the weak solution of the system with appropriate initial data.
\end{abstract}


\section{Introduction}
\setcounter{section}{1}\setcounter{equation}{0}

We consider in this paper the following 2D generalized quasi-geostrophic equation
\begin{equation}\label{gQG}
 \begin{cases}
  \partial_t \theta + u\cdot\nabla \theta +\nu \Lambda^{\beta}\theta =0, \\
  u= \Lambda^{\alpha}\mathcal{R}^{\perp}\theta,  \quad
  \theta(0,x)= \theta_0(x),
 \end{cases}
\end{equation}
where $\alpha\in ]0,1[$, $\beta\in ]0,2]$ and $\nu\geq 0$. The operator
$\Lambda^\beta:=(-\Delta)^{\frac{\beta}{2}}$ is defined through the Fourier transform
$$
 \widehat{\Lambda^\beta f}(\zeta) = |\zeta|^\beta \hat{f}(\zeta),
$$
and $\mathcal{R}^{\perp}:= (-\mathcal{R}_2,\mathcal{R}_1)$ are the usual Riesz transforms
\begin{equation*}
  \mathcal{R}_j\theta(x):=C\,\textrm{p.v.}\int_{\mathbb{R}^2}\frac{x_j-z_j}{|x-z|^{3}}\theta(z)\mathrm{d}z,\quad j=1,2.
\end{equation*}

\eqref{gQG} may be termed as the active scalar evolution equation in general. Here, we say the function $\theta(t,x)$ an active scalar if it satisfies
\begin{equation}\label{ActS0}
 \partial_t \theta + u\cdot\nabla \theta + \nu \Lambda^\beta \theta=0,\quad \theta(0,x)=\theta_0(x),
\end{equation}
where $\beta\in [0,2]$ and $u$ is defined by $\theta$ in some way. When $\alpha=0$ in \eqref{gQG}, it corresponds to the more well-known
active scalar equation------the quasi-geostrophic equation, which arises from the geostrophic study of the strongly rotating flows (cf. \cite{Con1}).
When $\alpha=1$ and $\beta=2$, although the flow term in \eqref{gQG} vanishes, we can still view another active scalar equation------the magnetogeostrophic equation introduced in \cite{FriV} as a meaningful generalization of this endpoint case,
where the divergence-free three-dimensional velocity $u$
satisfies $u= M [\theta]$ with $ M $ a nonlocal differential operator of order 1.
When $\alpha\in]-1,0[$, it is just the modified
quasi-geostrophic equation introduced in \cite{ConIW} by Constantin, Iyer and Wu.

\eqref{gQG} has the scaling invariance, i.e., if $\theta(t,x)$ is a solution of \eqref{gQG}, then
$$
\theta_\lambda(t,x)\rightarrow \lambda^{\beta-\alpha-1} \theta(\lambda^\beta t,\lambda x),\; \lambda>0
$$
is also a solution. Thus according to the $L^\infty$ maximum principle (cf. \cite{CorC}), we say $\beta>\alpha+1$, $\beta=\alpha+1$ and
$\beta<\alpha+1$ are the subcritical, critical and supercritical cases, respectively.
Note that when $\alpha=0$, it just corresponds to the classification of the classical quasi-geostrophic equation.

Compared with the quasi-geostrophic equation, \eqref{gQG} only has an additional operator $\Lambda^\alpha$ in $u$; but this positive derivative operator
will always produce much difficulty, so that many results can not (at least directly) extend to the generalized quasi-geostrophic equation \eqref{gQG}.
Before further proceeding, we first recall some noticeable results about quasi-geostrophic equation.
It has been known since \cite{Resnick} that the quasi-geostrophic equation has the global weak solutions for all cases $\beta\in ]0,2]$ and the global smooth
solution associated with suitable initial data for the subcritical case $\beta>1$.
See \cite{Con2,Wu 01} for other global results related to the subcritical case.
For the critical case $\beta=1$, Constantin et al in \cite{Con3} showed the global well-posedness
of the classical solution under the condition that the zero-dimensional $L^\infty$ norm of the data is small.
This smallness assumption was removed by Kiselev et al in \cite{KisNV}, where they
obtained the global well-posedness for the arbitrary periodic smooth data by using a new method that may be termed as
the nonlocal maximum principle method. Almost at the same time and from another totally different direction, Caffarelli and Vasseur \cite{Caffarelli}
resolved the problem to establish the global regularity of weak solutions
by deeply exploiting the De Giorgi method. See \cite{KisN} for a third but also quite different proof of the global regularity issue.
For the supercritical case $\beta<1$, although the problem that whether the equation has global regularity or not is still open, some partial
results have been proved. Local well-posedness for arbitrary initial data and global well-posedness under a smallness condition have been considered by many
authors (cf. \cite{ChenMZ,HmiK,Wu 04} and references therein).
From the direction of weak solutions, Constantin and Wu in \cite{ConW} showed a regularity criterion for the weak solutions in terms of the uniform (in time) H\"older estimates. Based on the criterion and considering the eventual regularity issue,
Silvestre in \cite{Sil} proved that for some $\beta<1/2$
but sufficiently close to $1/2$, the weak solutions become regular after a finite time; and then
Kiselev in \cite{Kis} developed the nonlocal maximum principle method
to show the eventual regularity for all the supercritical range $\beta\in ]0,1[$, see also \cite{Dab} for another proof by developing the method of \cite{KisN}.

The equation \eqref{gQG} with general $\alpha$ in $u$ was firstly introduced in \cite{ConIW}, where Constantin et al considered the critical case $\alpha\in ]-1,0[$ and $\beta=\alpha+1$ and they showed the global regularity of the weak solutions by applying Caffarelli-Vasseur's method.
Then the authors in \cite{MiaoXue} treated the
whole critical case $\alpha\in ]-1,1[$ and $\beta=\alpha+1$ by using the method of \cite{KisNV}: for the case $\alpha\in]-1,0[$,
global well-posedness of the smooth solution was proved; while for the case $\alpha\in ]0,1[$,
global well-posedness of the smooth solution was obtained under the condition of small $L^\infty$ initial data.
For the supercritical range $\alpha\in ]-1,0[$ and $\beta\in ]0,\alpha+1[$, Kiselev in \cite{Kis} also showed the eventual regularity of the weak solution.

In this article we focus on the regularity issue of the generalized quasi-geostrophic equation \eqref{gQG}.
Since \eqref{gQG} is more "singular" than the classical quasi-geostrophic equation,
we can not expect to obtain better regularity results than the corresponding results in the quasi-geostrophic equation.
In fact, we here prove the global regularity
for all the subcritical case, the local and eventual regularity for the case $\alpha\in ]0,1[$, $\beta\in ]2\alpha,\alpha+1]$.
First for the subcritical case, we precisely have
\begin{theorem}\label{thm1}
 Let $\nu>0$, $\alpha\in]0,1[$, $\beta\in]\alpha+1,2]$ and the initial data $\theta_{0}\in H^{m}(\mathbb{R}^2)$, $m>2$, then
 there exists a unique global solution
 \begin{equation*}
  \theta\in\mathcal{C}([0,\infty[;H^{m}) \cap L_{\mathrm{loc}}^2([0,\infty[; H^{m+\beta/2})\cap \mathcal{C}^{\infty}(]0,\infty[\times \mathbb{R}^{2})
 \end{equation*}
 to the generalized quasi-geostrophic equation \eqref{gQG}.
\end{theorem}

For the critical and supercritical cases, we need to introduce the notion of weak solution. Based on Theorem \ref{thm1}, it will be
convenient to consider the following system with additional viscosity
\begin{equation}\label{app-gQG}
 \begin{cases}
  \partial_t \theta + u\cdot\nabla \theta +\nu \Lambda^{\beta}\theta - \epsilon \Delta \theta =0, \\
  u= \Lambda^{\alpha}\mathcal{R}^{\perp}\theta.  \quad
  \theta(0,x)= \theta_0(x).
 \end{cases}
\end{equation}
Note that if the regularity of $\theta_0$ is very low (e.g. $\theta_0\in L^2$), we can replace $\theta_0(x)$ by the mollified data $\psi_\epsilon*\theta_0(x)$. Then passing to the limit $\epsilon\rightarrow 0$, we can get a weak solution for an appropriate range of $\alpha$, $\beta$ (see Proposition \ref{prop GlWeak} in the appendix). Then our main result is the following.

\begin{theorem}\label{thm2}
 Let $\nu>0$, $\alpha\in]0,1[$, $\beta\in ]2\alpha,\alpha+1]$ and $\theta_0\in H^m(\mathbb{R}^2)$ with $m>2$. Let $\theta(t,x)$ be the weak
 solution of \eqref{gQG} obtained by taking the limit $\epsilon\rightarrow 0$ of the corresponding smooth solution of \eqref{app-gQG}. Then
 there exist $0<T_1\leq T_2<\infty$ which depend only on $\nu,\alpha,\beta$ and $\theta_0$
 such that
 \begin{equation*}
   \theta \in \mathcal{C}^{\infty}(]0,T_1[\times \mathbb{R}^2)\cap \mathcal{C}^\infty (]T_2,\infty[\times\mathbb{R}^2).
 \end{equation*}
\end{theorem}

The local part of Theorem \ref{thm1} and Theorem \ref{thm2} corresponds to Proposition \ref{prop local} in the appendix,
and the global existence issue of the weak
solution is considered in Proposition \ref{prop GlWeak}. We use the newly developed method from \cite{KisNV,Kis} to treat the global part of Theorem \ref{thm1} and the eventual regularity part of Theorem \ref{thm2}. We first state the general criterion leading to the nonlocal maximum principles
for the whole space active scalars in the section \ref{sec GC}; then we apply the general criterion to prove the remaining parts of Theorem \ref{thm1}
and \ref{thm2} in the section \ref{sec glbol} and section \ref{sec EvtReg} respectively.

\begin{remark}\label{rem thm}
  In the global part of Theorem \ref{thm1}, we use the nonlocal maximum principle method to show the Lipschitz norm of the solution is
  uniformly bounded for all the existence time, then combining with the blowup criterion \eqref{eq blowup} leads to the global result.
  This is very different from the process
  in \cite{Con2,Wu 01}, and it seems hard to directly apply these classical methods to \eqref{gQG}.

  In the eventual regularity part of Theorem \ref{thm2},
  similarly as \cite{Kis,Dab,Sil}, we base on the regularity criterion Proposition \ref{prop RegC},
  and the regularization mechanism is different from the usual way
  to prove eventual regularity, which is by combining a global regularity result for small data with a suitable decay of some norm of the weak solution.
  Though in the proof of the critical case the
  decay (in time) of the $L^\infty$ norm is used, to obtain regularity we still need to wait
  an extra period of time after the $L^\infty$ norm is under control.
\end{remark}

\begin{remark}\label{rem thm1}
  Proposition \ref{prop RegC} indeed implies the regularity criterion of the weak solution for the critical and some supercritical equation \eqref{gQG} in terms of
  uniform (in time) H\"older estimates. We note that in the critical case of \eqref{gQG}, the criterion calls for that the regularity index $\sigma>\frac{\alpha}{2}$ with $\alpha\in ]0,1[$,
  which is essentially stronger than the criterion of the critical quasi-geostrophic equation (where $\sigma>0$). Thus if we rely on this criterion,
it will be not sufficient to obtain the global regularity
of the critical gQG equation \eqref{gQG} by applying the method of Caffarelli and Vasseur \cite{Caffarelli} (without introducing new ideas).
This is analogous with the case of applying the method of Kiselev et al \cite{KisNV} (cf. \cite{MiaoXue}). It is also interesting to point out
that a similar regularity criterion with the regularity index $\sigma>\frac{1}{2}$ should be considered for the critical model introduced in \cite{FriV}.
\end{remark}

\begin{remark}\label{rem thm2}
  The restriction $\beta>2\alpha$ is from Proposition \ref{prop local} and Proposition \ref{prop GlWeak} of constructing local smooth solution
  and global weak solution respectively. While in the proof of eventual regularity, under the condition
  $\beta>\alpha$ is already sufficient to show the uniform (in $\epsilon$) eventual $\mathcal{C}^{\tilde{\sigma}}$ regularity ($\tilde{\sigma}>1$)
  for the solution of \eqref{app-gQG}.
\end{remark}

\section{Preliminaries}
\setcounter{section}{2}\setcounter{equation}{0}

In this preparatory section, we present the definitions and some
related results of the Sobolev spaces, H\"older space and Besov spaces, also we
provide some important estimates which will be used later.

We begin with introducing some notations.
\\
$\diamond$ Throughout this paper, $C$ stands for a constant which may be different from line to line. We
sometimes use $A\lesssim B$ instead of $A\leq C B$, and use $A\lesssim_{\beta,\gamma\cdots}B$ instead of $A\leq C(\beta,\gamma,\cdots)B$ with $C(\beta,\gamma,\cdots)$ a constant depending on $\beta,\gamma,\cdots$. For $A\thickapprox B$ we mean $A\lesssim B\lesssim A$.
\\
$\diamond$ Denote by
$\mathcal{S}(\mathbb{R}^{n})$ the Schwartz space of rapidly
decreasing smooth functions, $\mathcal{S}'(\mathbb{R}^{n})$ the
space of tempered distributions,
$\mathcal{S}'(\mathbb{R}^{n})/\mathcal{P}(\mathbb{R}^{n})$ the
quotient space of tempered distributions which modulo polynomials.
\\
$\diamond$
$\mathcal{F}f$ or $\hat{f}$ denotes the Fourier transform, that is
$\mathcal{F}f(\zeta)=\hat{f}(\zeta)=\int_{\mathbb{R}^{n}}e^{-ix\cdot\zeta}f(x)\textrm{d} x,$
while $\mathcal{F}^{-1}f$ the inverse Fourier transform, namely,
$\mathcal{F}^{-1}f(x)=(2\pi)^{-n}\int_{\mathbb{R}^{n}}e^{ix\cdot\zeta}f(\zeta)\textrm{d}
\zeta$.
\\
$\diamond$ Denote by $B_r(x)$ the ball in $\mathbb{R}^n$ centered at $x$ with radius $r$. We abbreviate it by $B_r$, if the center is the origin.
Denote by $B^c_{r}(x)$ the complement set of $B_r(x)$ in $\mathbb{R}^n$.

Now we give the definition of ($L^2$-based) Sobolev space and H\"older space. For $s\in\mathbb{R}$, the inhomogeneous Sobolev space
\begin{equation*}
H^{s}:=\Big\{f\in \mathcal{S}'(\mathbb{R}^{n});
\norm{f}^{2}_{H^{s}}:=\int_{\mathbb{R}^{n}}(1+|\zeta|^{2})^{s}|\hat{f}(\zeta)|^{2}\textrm{d}
\zeta<\infty\Big\}
\end{equation*}
Also one can define the corresponding homogeneous space:
\begin{equation*}
\dot{H}^{s}:=\Big\{f\in
\mathcal{S}'(\mathbb{R}^{n})/\mathcal{P}(\mathbb{R}^{n});
\norm{f}^{2}_{\dot{H}^{s}}:=\int_{\mathbb{R}^{n}}|\zeta|^{2s}|\hat{f}(\zeta)|^{2}\textrm{d}
\zeta<\infty\Big\}
\end{equation*}
For $\delta\in]0,1[ $, the H\"older space $\mathcal{C}^\delta$ is the set of the tempered distribution $f$ such that
\begin{equation*}
  \|f\|_{\mathcal{C}^\delta}:= \|f\|_{L^\infty} + \sup_{x\neq y} \frac{|f(x)-f(y)|}{|x-y|^\delta}<\infty.
\end{equation*}

To define Besov space we need the following dyadic unity partition. Choose two nonnegative radial
functions $\chi$, $\varphi\in \mathcal{D}(\mathbb{R}^{n})$ be
supported respectively in the ball $\{\zeta\in
\mathbb{R}^{n}:|\zeta|\leq \frac{4}{3} \}$ and the shell $\{\zeta\in
\mathbb{R}^{n}: \frac{3}{4}\leq |\zeta|\leq
  \frac{8}{3} \}$ such that
\begin{equation*}
 \chi(\zeta)+\sum_{j\geq 0}\varphi(2^{-j}\zeta)=1, \quad
 \forall\zeta\in\mathbb{R}^{n}; \qquad
 \sum_{j\in \mathbb{Z}}\varphi(2^{-j}\zeta)=1, \quad \forall\zeta\neq 0.
\end{equation*}
For all $f\in\mathcal{S}'(\mathbb{R}^{n})$ we define the
nonhomogeneous Littlewood-Paley operators
\begin{equation*}
 \Delta_{-1}f:= \chi(D)f; \;\;
 \Delta_{j}f:= \varphi(2^{-j}D)f,\; S_j f:=\sum_{-1\leq k\leq j-1}\Delta_k f,\quad \forall j\in\mathbb{N},
\end{equation*}
And the homogeneous Littlewood-Paley operators can be defined as follows
\begin{equation*}
  \dot{\Delta}_{j}f:= \varphi(2^{-j}D)f;\; \dot S_j f:= \sum_{k\in\mathbb{Z},k\leq j-1}\dot \Delta_k f, \quad \forall j\in\mathbb{Z}.\quad
\end{equation*}

Now we introduce the definition of Besov spaces . Let $(p,r)\in
[1,\infty]^{2}$, $s\in\mathbb{R}$, the nonhomogeneous Besov space
\begin{equation*}
  B^{s}_{p,r}:=\Big\{f\in\mathcal{S}'(\mathbb{R}^{n});\norm{f}_{B^{s}_{p,r}}:=\norm{2^{js}\norm{\Delta
  _{j}f}_{L^{p}}}_{\ell^{r}}<\infty  \Big\}
\end{equation*}
and the homogeneous space
\begin{equation*}
  \dot{B}^{s}_{p,r}:=\Big\{f\in\mathcal{S}'(\mathbb{R}^{n})/\mathcal{P}(\mathbb{R}^{n});
  \norm{f}_{\dot{B}^{s}_{p,r}}:=\norm{2^{js}\norm{\dot{\Delta}
  _{j}f}_{L^{p}}}_{\ell^{r}(\mathbb{Z})}<\infty  \Big\}.
\end{equation*}
We point out that for all $s\in\mathbb{R}$, $B^{s}_{2,2}=H^{s}$ and
$\dot{B}^{s}_{2,2}=\dot{H}^{s}$.

The classical space-time Besov space $L^{\rho}([0,T],B^{s}_{p,r})$, abbreviated by
$L^{\rho}_{T}B^{s}_{p,r}$, is the set of tempered distribution $f$
such that
\begin{equation*}
  \norm{f}_{L^{\rho}_{T}B^{s}_{p,r}}:=\norm{\norm{2^{js}\norm{\Delta_{j}f}_{L^{p}}}_{\ell^{r}}}_{L^{\rho}([0,T])}<\infty.
\end{equation*}
We can similarly extend to the homogeneous one $L^{\rho}_{T}\dot{B}^{s}_{p,r}$.

Bernstein's inequality is fundamental in the analysis involving
Besov spaces (see \cite{Lemarie})
\begin{lemma}
Let $f\in L^{a}$, $1\leq a\leq b\leq \infty$. Then for every $
(k,q)\in\mathbb{N}^{2}$ there exists a constant $C>0$ such that
\begin{equation*}
 \sup_{|\alpha|=k}\norm{\partial^{\alpha} S_{q}f}_{L^{b}}\leq C 2
 ^{q(k+n(\frac{1}{a}-\frac{1}{b}))}\norm{f}_{L^{a}},
\end{equation*}
\begin{equation*}
 C^{-1}2
 ^{q k}\norm{f}_{L^{a}}\leq \sup_{|\alpha|=k}\norm{\partial^{\alpha}
 \Delta_{q}f}_{L^{a}}\leq C 2^{qk}\norm{f}_{L^{a}}
 \end{equation*}
\end{lemma}

Next we state an important maximum principle for the transport-diffusion equation (cf. \cite{CorC}).
\begin{proposition}\label{prop MP}
Let $u$ be a smooth divergence-free vector field and $f$ be a smooth
function. Assume that $\theta$ is the smooth solution of the following equation
\begin{equation*}
  \partial_{t}\theta + u\cdot\nabla\theta+\nu \Lambda^{\beta}\theta=f,\quad \mathrm{div}u =0,
\end{equation*}
with initial datum $\theta_0$ and $\nu\geq 0$, $0\leq\beta\leq 2$, then for every $p\in
[1,\infty]$ we have
\begin{equation}\label{eq TDMaxPrin}
 \norm{\theta(t)}_{L^{p}}\leq
 \norm{\theta_{0}}_{L^{p}}+\int^{t}_{0}\norm{f(\tau)}_{L^{p}}
 \,\textrm{d}\tau.
\end{equation}
\end{proposition}

Finally we concern a uniform decay estimate of the global smooth solution of \eqref{app-gQG}.

\begin{lemma}\label{lem LinfDC}
Let $\alpha\in ]0,1[$, $\beta\in ]0,2]$ and $\theta(t,x)$ be the global smooth solution of the equation \eqref{app-gQG}
\begin{equation*}
  \partial_t \theta + u\cdot\nabla \theta +\nu \Lambda^{\beta}\theta - \epsilon \Delta \theta =0, \;\;
  u= \Lambda^{\alpha}\mathcal{R}^{\perp}\theta, \quad \theta(0,x)=\theta_0(x),
\end{equation*}
with initial data $\theta_0(x)\in H^m$, $m>2$. Then we have the decay estimate
\begin{equation*}
  \|\theta(t)\|_{L^\infty} \leq \frac{\|\theta_0\|_{L^\infty}}{\Big(1+ C\frac{\|\theta_0\|_{L^\infty}^\beta}{\|\theta_0\|_{L^2}^\beta}t\Big)^{1/\beta}},\quad t>0,
\end{equation*}
where $C$ is an absolute constant independent of the value of $\epsilon$.
\end{lemma}

The proof is quite similar to the proof of Theorem 4.1 and Corollary 4.2 in \cite{CorC}, thus we omit it.

\section{The Modulus of Continuity and the General Criterion}\label{sec GC}
\setcounter{section}{3}\setcounter{equation}{0}

Motivated by \cite{Kis}, in this section we state a general criterion leading to the nonlocal maximum principles for the whole space active scalars.
Here we call the function $\theta(t,x)$ whole space active scalar if the space variable of the active scalar $\theta(t,x)$ is over the whole space $\mathbb{R}^n$.

We begin with introducing some terminology (cf. \cite{Kis}).
\begin{definition}
  A function $\omega:]0,\infty[ \rightarrow ]0,\infty[$ is called a modulus of continuity (MOC) if $\omega$ is continuous on $]0,\infty[$, increasing, concave, and
  piecewise $\mathcal{C}^2$ with one-sided derivatives defined at every point in $]0,\infty[$ (maybe infinite at $\xi=0$). We call that a function $f:\mathbb{R}^n\rightarrow \mathbb{R}^l$ obeys the modulus of continuity $\omega$ if $|f(x)-f(y)| < \omega(|x-y|)$ for every $x\neq y\in \mathbb{R}^n$.
\end{definition}

Notice that from the definition, the inverse function of $\omega(\xi)$ is uniquely determined, and denote it by $\omega^{-1}(y)$. Clearly $\omega(\omega^{-1}(y))=y$ for all $\omega^{-1}(y)<\infty$, $\omega^{-1}(\omega(\xi))=\xi$, and $\omega^{-1}(y)$ is continuous, increasing and convex. If we consider the time-dependent MOC $\omega(t,\xi)$, we shall correspondingly denote $\omega^{-1}(t,y)$ as its inverse function.

Then the main result in this section is as follows.
\begin{proposition}\label{prop GC}
 Let $\theta(t,x)$ be a smooth solution of the following whole space active scalar equation
\begin{equation}\label{appActS}
 \partial_t \theta + u\cdot\nabla \theta + \nu \Lambda^\beta \theta -\epsilon \Delta \theta=0,\quad \theta(0,x)=\theta_0(x), \; x\in\mathbb{R}^n,
\end{equation}
 with $\epsilon\geq 0$, $\nu>0$, $\beta\in]0,2]$ and $\theta$ have a suitable spacial decay property (i.e. $\|\nabla \theta(t)  \|_{L^\infty(B^c_{R})}\rightarrow 0$ as $R\rightarrow \infty$ for every $t>0$). Assume that for every $t\geq 0$, $\omega(t,\xi)$ is a modulus of continuity and satisfies $\omega^{-1}(t,\frac{2}{\delta_1}\|\theta_0 \|_{L^\infty})<\infty$ with some fixed number $\delta_1\in [\frac{2}{3},1[$; and that for every fixed point $\xi$, $\omega(t,\xi)$ is piecewise $\mathcal{C}^1$ in time variable with one-sided derivatives defined at all points, and in particular that for all $\xi$ near infinity, $\omega(t,\xi)$ is continuous in $t$ uniformly in $\xi$ (i.e. $\exists M>0$ a large number, such that $|\omega(t+h,\xi)-\omega(t,\xi)|\leq \tilde{\epsilon}$ for all $\xi\geq M$ with $h$ depending only on $M,t,\tilde{\epsilon}$). Assume in addition that $\omega(t,0+)$ and $\partial_\xi \omega(t,0+)$ is continuous in $t$ with values in $\mathbb{R}\cup \infty$, and that one of the three conditions below is satisfied
\\
(a) for every $t\geq 0$, $\omega(t,0+)>0$,
\\
(b) for every $t\geq 0$, $\omega(t,0+)=0$, $\partial_\xi\omega( t,0+)=\infty$,
\\
(c) for every $t\geq 0$, $\omega(t,0+)=0$, $\partial_\xi\omega( t,0+)<\infty$, $\partial_{\xi\xi}\omega (t,0+) =-\infty$.
\\
Let the initial data $\theta_0(x)$ obey $\omega(0,\xi):=\omega_0(\xi)$.

Then $\theta(T,x)$ obeys the modulus of continuity $\omega(T,\xi)$ provided that $\omega(t,\xi)$ satisfies
\begin{equation}\label{eq keyGC}
  \partial_t \omega(t,\xi)> \Omega(t,\xi)\partial_\xi\omega(t,\xi) + \nu \Upsilon_\beta(t,\xi) + 2\epsilon \partial_{\xi\xi}\omega(t,\xi),
\end{equation}
for all $0<t\leq T$ and $\xi>0$ such that $\omega(t,\xi)\leq 2\| \theta(t,\cdot)\|_{L^\infty}$, and
where $\Omega(t,\xi)$ is from the bound that for every $x\in\mathbb{R}^n$ and every unit vector $e\in \mathbb{R}^n$,
\begin{equation}\label{eq GC-Ome}
  |(u(t,x+\xi e)-u(t,x))\cdot e|\leq \Omega(t,\xi),
\end{equation}
and $\Upsilon_\beta(t,\xi)$ is (usually) given by
\begin{equation}\label{eq GC-Ups}
\begin{split}
  \Upsilon_\beta(t,\xi) =&  c_\beta
 \int_{0}^{\frac{\xi}{2}}\frac{\omega(t,\xi+2\eta)+\omega (t,\xi-2\eta)
 -2\omega (t,\xi)}{\eta^{1+\beta}}\textrm{d} \eta \\
    & +c_\beta   \int_{\frac{\xi}{2}}^{\infty}\frac{\omega (t,2\eta+\xi)-\omega (t,2\eta-\xi)
 -2\omega (t,\xi)}{\eta^{1+\beta}}\textrm{d}\eta,
\end{split}
\end{equation}
with $c_\beta$ an absolute constant depending only on $\beta$ and $n$.
In \eqref{eq keyGC}, at the points where $\partial_t\omega(t,\xi)$ ($\partial_\xi\omega(t,\xi)$) does not exist, the smaller (larger) value of the one-sided derivative should be taken.
\end{proposition}


\begin{proof}
 We proceed by contradiction. Suppose that $\theta(t,x)$ no longer obeys $\omega(t,\xi)$ for some $t>0$, then we claim that
 there must exist $T_*>0$ and two fixed points $\bar{x}\neq \bar{y}\in \mathbb{R}^n$ such that $\theta(t,x)$ obeys $\omega(t,\xi)$ for every $t<T_*$, while
 \begin{equation}\label{eq scenario}
   \theta(T_*,\bar{x})-\theta(T_*,\bar{y})=\omega(T_*,|\bar{x}-\bar{y}|).
 \end{equation}
 Indeed, denote
 \begin{equation}
   T_*:=\sup\{T\in[0,\infty[; \forall t <T, \, \theta(t,x)\, \textrm{obeys} \, \omega(t,\xi)  \}
 \end{equation}
 in other words, $T_*$ is the minimal time that $\theta(t,x)$ no longer obeys $\omega(t,\xi)$. Clearly, we see that
 $\theta(T_*,x)-\theta(T_*,y)\leq \omega(T_*,|x-y|)$ for all $x,y\in\mathbb{R}^n$, otherwise from the time
 continuity property we shall get a contradiction (cf. \cite{Kis}). Then for all $x\neq y\in\mathbb{R}^n$, set
 $$F(t,x,y):=\frac{|\theta(t,x)-\theta(t,y)|}{\omega(t,|x-y|)}.$$
 Clearly $F(T_*,x,y)\leq 1$, and we assume that $F(T_*,x,y)< 1$ for all $x\neq y$, otherwise the claim is proved. We shall prove that there exists some small
 $h>0$ such that $F(t,x,y)<1$ for all $x,y\in\mathbb{R}^n,x\neq y$ and $t\in[T_*,T_*+h]$, which contradicts with the definition of $T_*$.

 First, denote $C_0(t):=\omega^{-1}(t,\frac{2}{\delta_1}\norm{\theta_0}_{L^\infty})$ with some $\delta_1\in [2/3,1[$. When
 $|x-y|\geq C_0(T_*)$, we have
 \begin{equation*}
    2\norm{\theta_0}_{L^\infty} =\delta_1 \omega(T_*,C_0(T_*))\leq \delta_1 \omega(T_*,|x-y|).
 \end{equation*}
 From the uniform continuity property of $\omega(t,\xi)$, there exist a small $h_1>0$ depending only on $M,T_*,C_0(T_*),\delta_1$ such that
 for all $|x-y|\geq C_0(T_*)$ and $t\in [T_*,T_*+h_1]$
 \begin{equation}
   |\theta(t,x)-\theta(t,y)|\leq 2\norm{\theta_0}_{L^\infty} \leq \frac{\delta_1+1}{2} \omega(t,|x-y|).
 \end{equation}
 Since $\theta$ is smooth and has a spacial decay property, then for every $\tilde{\epsilon}>0$, there exist $h_2,R(T_*)>0$ depending on $\tilde{\epsilon},T_*$ such that for every $t\in [T_*,T_*+h_2]$ we have
 \begin{equation*}
   \norm{\nabla \theta(t)}_{L^\infty(B^c_{R(T_*)})} \leq \norm{\nabla \theta(T_*)}_{L^\infty(B^c_{R(T_*)})} +\frac{\tilde{\epsilon}}{2} \leq \tilde{\epsilon}.
 \end{equation*}
 Hence for all $|x-y|\leq C_0(T_*)$ and $x$ or $y$ belongs $B^c_{R(T_*)+C_0(T_*)}$, we have for every $t\in [T_*,T_*+h_2]$
 \begin{equation*}
   |\theta(t,x)-\theta(t,y)|\leq \norm{\nabla \theta(t)}_{L^\infty(B^c_{R(T_*)})} |x-y| \leq \tilde{\epsilon} |x-y|.
 \end{equation*}
 Note that from the concavity of $\omega(t,\cdot)$, we get $\frac{\omega(t,C_0(T_*))}{C_0(T_*)}|x-y|\leq \omega(t, |x-y|)$, thus it only
 needs to choose $\tilde{\epsilon}\leq \frac{1}{2} \frac{\omega(t,C_0(T_*))}{C_0(T_*)}$ so that
 \begin{equation}
   |\theta(t,x)-\theta(t,y)|\leq (1/2)\omega(t,|x-y|).
 \end{equation}
 Next it suffices to consider $x,y\in B_{R(T_*)+C_0(T_*)}$. In a similar way as treating the corresponding part in \cite{Kis}, we get that
 there exist $\kappa,h_3>0$, $\rho\in ]0,1[$ such that for all $x, y\in B_{R(T_*)+C_0(T_*)}$, $x\neq y$, $|x-y|\leq \kappa$ and $t\in [T_*,T_*+h_3]$
 \begin{equation*}
   F(t,x,y)<\rho.
 \end{equation*}
 Then it remains to consider the continuous function $F(t,x,y)$ on the compact set
 $$ \mathcal{K}:= \{(x,y)\in\mathbb{R}^n\times\mathbb{R}^n; \max\{|x|,|y|\}\leq R(T_*)+C_0(T_*),\,|x-y|\geq \kappa \}.$$
 From $F(T_*,x,y)<1$ on $\mathcal{K}$ and the continuity in time of $F$, there exists a small number $h_4>0$ such that the strict inequality holds
 for every $t\in [T_*,T_*+h_4]$ and $(x,y)\in\mathcal{K}$. Let $h:=\min\{h_1,h_2,h_3,h_4\}$, then $F(t,x,y)<1$ for all $t\in[T_*,T_*+h]$
 and $x,y\in\mathbb{R}^n,x\neq y$, contradicting the definition of $T_*$.
 Therefore, there have to exist two fixed points $\bar{x},\bar{y}$ (in fact in $\mathcal{K}$) such that \eqref{eq scenario} holds.

 Now in this breakdown scenario, we shall use the equation and calculate the time derivative of $\theta(t,\bar{x})-\theta(t,\bar{y})$ at time $t=T_*$. Indeed, since $\theta$ is smooth, from \eqref{appActS} we have
 \begin{equation*}
 \begin{split}
   \partial_t(\theta(t,\bar{x})-\theta(t,\bar{x}))|_{t=T_*} = &
   -(u\cdot\nabla)\theta(T_*,\bar{x})+ (u\cdot\nabla)\theta(T_*,\bar{y})-\\ & -\nu\Lambda^\beta \theta(T_*,\bar{x})+\nu \Lambda^\beta \theta(T_*,\bar{y})
   +\epsilon \Delta\theta(T_*,\bar{x})-\epsilon\Delta\theta(T_*,\bar{y})
 \end{split}
 \end{equation*}
 Denote $\bar{\xi}=|\bar{x}-\bar{y}|$, $\ell=\frac{\bar{x}-\bar{y}}{|\bar{x}-\bar{y}|}$ and $v$ the arbitrary unit vector orthogonal to $\ell$. Then almost parallel to
 Proposition 2.4 in \cite{Kis} we get
 \begin{equation*}
   \partial_\ell\theta(T_*,\bar{x})=\partial_\ell \theta(T_*,\bar{y})=\partial_\xi\omega(\bar{\xi}),\quad \partial_{\ell\ell} \theta(T_*,\bar{x}) \leq\partial_{\xi\xi}\omega(T_*,\bar{\xi}),\; \partial_{\ell\ell}\theta(T_*,\bar{y})\geq -\partial_{\xi\xi}\omega(T_*,\bar{\xi});
 \end{equation*}
 \begin{equation*}
   \partial_v \theta(T_*,\bar{x})=\partial_v \theta(T_*,\bar{y})=0,\quad \partial_{v v}\theta(T_*,\bar{x})-\partial_{v v}\theta(T_*,\bar{y})\leq 0.
 \end{equation*}
So we have
 \begin{equation*}
 \begin{split}
   |(u\cdot\nabla)\theta(T_*,\bar{x})-(u\cdot\nabla)\theta(T_*,\bar{y})|= |(u(T_*,\bar{x})-u(T_*,\bar{y}))\cdot \ell|\partial_\xi\omega(T_*,\bar{\xi})\leq
   \Omega(T_*, \bar{\xi})\partial_\xi\omega(T_*,\bar{\xi});
 \end{split}
 \end{equation*}
 and
 \begin{equation}\label{eq DisEpsilon}
   \epsilon \Delta\theta(T_*,\bar{x})-\epsilon\Delta\theta(T_*,\bar{y})\leq 2\epsilon\partial_{\xi\xi}\omega(T_*,\bar{\xi}).
 \end{equation}
 And the contribution of the dissipative term can be estimated as (cf. \cite{MiaoXue,KisNV})
 \begin{equation}\label{eq Upsi}
    -\nu\Lambda^\beta\theta(T_*,\bar{x})+\nu \Lambda^\beta\theta(T_*,\bar{y})\leq \nu \Upsilon_\beta(T_*,\bar{\xi});
 \end{equation}
 Note that in the proof of \eqref{eq Upsi}, there is not much difference between the case $\beta=2$ and the other case $\beta\in ]0,2[$, only observing that
 $\lim_{h\rightarrow 0}\frac{1}{h} G_h(\eta)= \lim_{h\rightarrow0}\frac{1}{h} \frac{1}{(4\pi h)^{1/2}}e^{-\eta^2/(4h)}\leq C \frac{1}{\eta^3}$;
 thus \eqref{eq Upsi} also offers another estimation of the LHS of \eqref{eq DisEpsilon}.

 Hence, based on the above analysis, we have
 \begin{equation*}
   \partial_t \Big( \frac{\theta(t,\bar{x})-\theta(t,\bar{y})}{\omega(t,\bar{\xi})}\Big)\Big|_{t=T_*}\leq
   \frac{\Omega(T_*,\bar{\xi})\partial_\xi \omega(T_*,\bar{\xi}) + \nu\Upsilon_\beta(T_*,\bar{\xi})+2\epsilon\partial_{\xi\xi}\omega(T_*,\bar{\xi})-\partial_t\omega(T_*,\bar{\xi})}{\omega(T_*,\bar{\xi})}.
 \end{equation*}
 From \eqref{eq keyGC}, the RHS of the above inequality is strictly negative, which clearly yields a contradiction with the choice of $T_*$.

\end{proof}

\section{Global well-posedness for $\alpha\in ]0,1[$ and $\beta\in ]\alpha+1,2]$}\label{sec glbol}
\setcounter{section}{4}\setcounter{equation}{0}

From Proposition \ref{prop local} in the Appendix, we assume that $T^*$ is the maximal existence time of the solution of \eqref{gQG} in
$\mathcal{C}([0,T^*[,H^m)\cap L^2([0,T^*[,H^{m+\frac{\beta}{2}})$. We shall apply the general criterion Proposition \ref{prop GC} to
show that some appropriate moduli of continuity are persisted, which implies that the Lipschitz
norm of the solution is bounded uniformly in time. Of course, this combined with \eqref{eq blowup} further leads to $T^*=\infty$.

In fact, from the scaling transformation of \eqref{gQG}, we shall find some stationary MOC
\begin{equation}\label{eq omelmd}
  \omega_\lambda(\xi):=\lambda^{\beta-\alpha-1}\omega(\lambda \xi),\quad \lambda\in]0,\infty[
\end{equation}
satisfying condition (c) and $\omega_{\lambda}^{-1}(3\|\theta_0\|_{L^\infty})<\infty$.

With no loss of generality, we assume that there is a fixed constant $c_0>0$ such that $\lim_{\xi\rightarrow \infty}\omega(\xi)>c_0$, that is $\omega^{-1}(c_0)<\infty$.
Then we can choose some $\lambda\in ]0,\infty[$ such that
\begin{equation}\label{eq lamdCod1}
\lambda^{\beta-\alpha-1}>\frac{3\|\theta_0\|_{L^\infty}}{c_0},
\end{equation}
and from $\omega_\lambda^{-1}(y)=\frac{1}{\lambda}\omega^{-1}(\frac{y}{\lambda^{\beta-\alpha-1}})$, we get $\omega^{-1}_{\lambda}(3\|\theta_0\|_{L^\infty})<\infty$. Now we check the condition that $\theta_0(x)$ obeys $\omega_\lambda(\xi)$ for appropriate $\lambda$.
First, from \eqref{eq lamdCod1}, we know that for every $x,y$ such that $\lambda |x-y| \geq \omega^{-1}(c_0)$
\begin{equation}
  |\theta_0(x)-\theta_0(y)|\leq 2\|\theta_0\|_{L^\infty}\leq \frac{2}{3} \lambda^{\beta-\alpha-1} c_0\leq \frac{2}{3}\omega_\lambda(|x-y|).
\end{equation}
Second, using the mean value theorem, we have
\begin{equation*}
  |\theta_0(x)-\theta_0(y)|\leq \|\nabla\theta_0\|_{L^\infty}|x-y|.
\end{equation*}
Let $0<\delta_0<\omega^{-1}(c_0)$. Due to the concavity of $\omega$, we infer that for every $x,y$ such that $\lambda|x-y|\leq \delta_0$,
\begin{equation*}
  \frac{\lambda^{\beta-\alpha-1}\omega(\delta_0)}{\delta_0}\leq \frac{\omega_\lambda(|x-y|)}{\lambda |x-y|}.
\end{equation*}
Thus by choosing $\lambda$ such that
\begin{equation}\label{eq lamdCod2}
  \lambda^{\beta-\alpha}> \frac{\delta_0}{\omega(\delta_0)} \|\nabla\theta_0\|_{L^\infty},
\end{equation}
we get that for every $x,y$ satisfying $x\neq y$ and $\lambda|x-y|\leq \delta_0$,
\begin{equation}
  |\theta_0(x)-\theta_0(y)|<\omega_\lambda(|x-y|).
\end{equation}
Finally, we consider the case $\delta_0\leq \lambda|x-y|\leq \omega^{-1}(c_0)$. Notice that $|\theta_0(x)-\theta_0(y)|\leq \frac{\omega^{-1}(c_0)}{\lambda}\|\nabla\theta_0\|_{L^\infty}$ and $\lambda^{\beta-\alpha-1}\omega(\delta_0)\leq \omega_\lambda(|x-y|)$.
Thus by choosing $\lambda$ satisfying
\begin{equation}\label{eq lamdCod3}
  \lambda^{\beta-\alpha}> \frac{\omega^{-1}(c_0)}{\omega(\delta_0)} \|\nabla\theta_0\|_{L^\infty},
\end{equation}
we obtain that for every $x,y$ satisfying $\delta_0\leq \lambda|x-y|\leq \omega^{-1}(c_0)$,
\begin{equation}
  |\theta_0(x)-\theta_0(y)|<\omega_\lambda(|x-y|).
\end{equation}
Hence, to fit our purpose, we can pick
\begin{equation}\label{eq lambda}
 \lambda= \max\Big\{\big(\frac{4\|\theta_0\|_{L^\infty}}{c_0}\big)^{\frac{1}{\beta-\alpha-1}}, \frac{\omega^{-1}(c_0)}{\|\theta_0\|_{L^\infty}} \|\nabla\theta_0\|_{L^\infty} \Big\},
\end{equation}
and $\delta_0=\omega^{-1}(\frac{2\|\theta_0\|_{L^\infty}}{\lambda^{\beta-\alpha-1}})$.

Then it remains to check \eqref{eq keyGC} for such $\omega_\lambda(\xi)$ with $\lambda$ given by \eqref{eq lambda}. For the contribution of the nonlinear term, from Lemma 3.2 in \cite{MiaoXue},
we know: if $\theta(t,x)$ obeys $\omega(\xi)$, then $u= \Lambda^\alpha \mathcal{R}^\bot \theta$ satisfies that for every $x,y\in\mathbb{R}^2$
\begin{equation}\label{eq uEsRu}
  |u(t,x)-u(t,y)|\leq \Omega_1(\xi),\quad \textrm{with} \quad  \xi:=|x-y|,
\end{equation}
where
\begin{equation}
   \Omega_1(\xi):= c_\alpha \Big( \int_0^\xi \frac{\omega(\eta)}{\eta^{1+\alpha}}\mathrm{d}\eta + \xi\int_\xi^\infty
   \frac{\omega(\eta)}{\eta^{2+\alpha}}\mathrm{d}\eta\Big)
\end{equation}
with $c_\alpha$ an absolute constant. Thus for such $\omega_\lambda(\xi)$ given by \eqref{eq omelmd}, correspondingly, by changing of variable we get
\begin{equation}\label{eq Ome1lmd}
  \Omega_{1,\lambda}(\xi)= c_\alpha \Big( \int_0^\xi \frac{\omega_\lambda(\eta)}{\eta^{1+\alpha}}\mathrm{d}\eta + \xi\int_\xi^\infty
   \frac{\omega_\lambda(\eta)}{\eta^{2+\alpha}}\mathrm{d}\eta\Big)= \lambda^{\beta-1}\Omega_1(\lambda\xi).
\end{equation}
While for the dissipative term, from \eqref{eq GC-Ups}, we see
\begin{equation}
\begin{split}
  \Upsilon_\beta(\xi) =   c_\beta
 \Big(\int_{0}^{\frac{\xi}{2}}\frac{\omega(\xi+2\eta)+\omega (\xi-2\eta)
 -2\omega (\xi)}{\eta^{1+\beta}}\textrm{d} \eta +
       \int_{\frac{\xi}{2}}^{\infty}\frac{\omega (2\eta+\xi)-\omega (2\eta-\xi)
 -2\omega (\xi)}{\eta^{1+\beta}}\textrm{d}\eta\Big),
\end{split}
\end{equation}
thus for $\omega_\lambda(\xi)$, we have
\begin{equation*}
\Upsilon_{\beta,\lambda}(\xi)=\lambda^{2\beta-\alpha-1}\Upsilon_\beta(\lambda\xi).
\end{equation*}
Then \eqref{eq keyGC} reduces to
\begin{equation*}
  \Omega_{1,\lambda}(\xi)\omega'_\lambda(\xi)+ \nu\Upsilon_{\beta,\lambda}(\xi)<0, \quad \textrm{for all}\; \xi>0,
\end{equation*}
equivalently,
\begin{equation*}
  \lambda^{2\beta-\alpha-1}(\Omega_1\omega'+ \nu\Upsilon_{\beta})(\lambda\xi)<0, \quad \textrm{for all}\; \xi>0.
\end{equation*}

Next, we shall construct a suitable modulus of continuity in the spirit of \cite{KisNV}.
 Choose two small positive numbers
$0<\gamma<\delta<1 $ and define the continuous functions $\omega$ as follows
\begin{equation}\label{eq modulus}
\begin{cases}
 \omega(\xi)=\xi-\xi^{\frac{3}{2}}
 \quad & \text{if} \quad 0< \xi\leq \delta, \\
 \omega'(\xi)=\frac{  \gamma}{4(\xi+\xi^{\beta})}
 \quad & \text{if} \quad \xi>\delta,
\end{cases}
\end{equation}
Note that, for small $\delta$, the left derivative of $\omega$ at
$\delta$ is about 1, while the right derivative equals
$\frac{\gamma}{4(\delta+\delta^{\beta})}<\frac{1}{4}$.
So $\omega$ is concave if $\delta$ is small enough (e.g. $\delta\leq \frac{1}{9}$). Clearly,
$\omega(0+)=0$, $\omega'(0+)=1$ and $\omega''(0+)= -\infty$. Due to $\beta>1$, $\omega$ is a bounded function, and at least $\lim_{\xi\rightarrow \infty}\omega(\xi)>\omega(\delta)= \delta-\delta^{\frac{3}{2}}$.

Then our target is to show that, for this MOC $\omega$,
\begin{equation*}
  \Omega_1(\xi)\omega'(\xi)+ \nu \Upsilon_\beta(\xi)<0  \quad \textrm{for all}\quad\xi>0.
\end{equation*}
More precisely, we need to prove the inequality
\begin{equation*}
\begin{split}
 c_\alpha \bigg[\int^{\xi}_{0}\frac{\omega(\eta)}{\eta^{1+\alpha}}\textrm{d} \eta+
 \xi\int_{\xi}^{\infty}\frac{\omega(\eta)}{\eta^{2+\alpha}}\textrm{d}
 \eta\bigg]\omega'(\xi)
 + c_\beta\nu \int_{0}^{\frac{\xi}{2}}\frac{\omega(\xi+2\eta)+\omega(\xi-2\eta)-2\omega(\xi)}{\eta^{1+\beta}}\textrm{d}
 \eta &\\
 +   c_\beta \nu \int_{\frac{\xi}{2}}^{\infty}\frac{\omega(2\eta+\xi)-\omega(2\eta-\xi)-2\omega(\xi)}{\eta^{1+\beta}}\textrm{d}
 \eta <0 \quad \text{for all}\quad \xi>0&.
\end{split}
\end{equation*}
\vskip0.2cm
Case 1: $0<\xi\leq \delta$.
\vskip0.2cm
Since $\frac{\omega(\eta)}{\eta}\leq \omega'(0+)=1$ for all $\eta>0$, we have
$$\int_{0}^{\xi}\frac{\omega(\eta)}{\eta^{1+\alpha}}\textrm{d} \eta \leq
\int_{0}^{\xi}\frac{1}{\eta^{\alpha}}\textrm{d} \eta \leq \frac{1}{1-\alpha} \xi^{1-\alpha},$$
 and
$$\int_{\xi}^{\delta}\frac{\omega(\eta)}{\eta^{2+\alpha}}\textrm{d}
 \eta \leq \int_{\xi}^{\delta}\frac{1}{\eta^{1+\alpha}}\textrm{d}
 \eta \leq \frac{1}{\alpha} \xi^{-\alpha}.$$
  Further,
\begin{equation*}
\begin{split}
 \int_{\delta}^{\infty}\frac{\omega(\eta)}{\eta^{2+\alpha}}\textrm{d}
 \eta & =\frac{1}{\alpha+1}\frac{\omega(\delta)}{\delta^{\alpha+1}}+
 \frac{1}{\alpha+1}\int_{\delta}^{\infty}\frac{\gamma}{4\eta^{\alpha+1}(\eta+\eta^{\beta})}\textrm{d}
 \eta \\
 & \leq \frac{1}{\alpha+1} \frac{1}{\delta^{\alpha}}+
  \frac{ \gamma}{4(\alpha+1)^2}\frac{1}{\delta^{\alpha+1}}\leq \frac{2 }{\delta^{\alpha}}\leq  2 \xi^{-\alpha}.
\end{split}
\end{equation*}
Obviously $\omega'(\xi)\leq\omega'(0)=1$, so we get that the positive
part is bounded by $ c_\alpha \xi^{1-\alpha}\frac{2}{\alpha(1-\alpha)}$.

For the negative part, we have
\begin{equation*}
\begin{split}
 c_\beta \int_{0}^{\frac{\xi}{2}}\frac{\omega(\xi+2\eta)+\omega(\xi-2\eta)-2\omega(\xi)}{\eta^{1+\beta}}\textrm{d}
 \eta \leq
 c_\beta\int_{0}^{\frac{\xi}{2}}\frac{\omega''(\xi)2\eta^{2}}{\eta^{1+\beta}}\textrm{d}
 \eta.
\end{split}
\end{equation*}
Due to that $\omega''(\xi)=-\frac{3}{4}\xi^{-\frac{1}{2}}<0$, we infer that the last expression is bounded by
\begin{equation*}
  -\frac{3}{2}c_\beta\xi^{-\frac{1}{2}} \int_{\frac{\xi}{4}}^{\frac{\xi}{2}}\frac{1}{\eta^{\beta-1}} \mathrm{d}\eta\leq
  -\frac{3}{2}c_\beta\xi^{-\frac{1}{2}} \frac{\xi}{4}(\frac{\xi}{2})^{1-\beta}\leq -\frac{3}{8}c_\beta \xi^{\frac{3}{2}-\beta}.
\end{equation*}

But from $\beta>\alpha+1$, clearly $\xi^{1-\alpha}\Big( \frac{2 c_\alpha}{\alpha(1-\alpha)} -\frac{3\nu c_\beta}{8} \xi^{\frac{1}{2}+\alpha-\beta}\Big)<0$ on
 $(0,\delta]$ when $\delta$ is small enough \big(i.e. $\delta^{\beta-\alpha-1/2}< \frac{3\nu\alpha(1-\alpha)c_\beta}{16 c_\alpha}$\big).
\vskip0.2cm
 Case 2: $\xi\geq \delta$
\vskip0.2cm

For $0\leq\eta\leq\delta$ we still have
$\omega(\eta)\leq \eta$ and for $\delta\leq
\eta\leq\xi$ we have $\omega(\eta)\leq\omega(\xi)$, then
\begin{equation*}
 \int_{0}^{\xi}\frac{\omega(\eta)}{\eta^{\alpha+1}}\textrm{d} \eta\leq
 \frac{\delta^{1-\alpha}}{1-\alpha}+\frac{\omega(\xi)}{\alpha}\Big(\delta^{-\alpha}-\xi^{-\alpha}\Big)\leq
 \frac{2\delta^{-\alpha}}{\alpha(1-\alpha)}\omega(\xi),
\end{equation*}
where the last inequality is due to that
$\frac{\delta}{2}\leq \delta-\delta^{\frac{3}{2}}=\omega(\delta)\leq\omega(\xi)$ for $\delta<\frac{1}{4}$.
Also
\begin{equation*}
\begin{split}
 \int_{\xi}^{\infty}\frac{\omega(\eta)}{\eta^{2+\alpha}}\textrm{d}
 \eta & =\frac{1}{\alpha+1}\frac{\omega(\xi)}{\xi^{\alpha+1}}+
 \frac{1}{\alpha+1}\int_{\xi}^{\infty}\frac{ \gamma}{4\eta^{\alpha+1}(\eta+\eta^{\beta})}\textrm{d} \eta \\
 & \leq \frac{1}{\alpha+1} \frac{\omega(\xi)}{\xi^{\alpha+1}}+   \frac{ \gamma}{4(\alpha+1)^2}\frac{1}{\xi^{\alpha+1}}\leq
 \frac{2\omega(\xi)}{\xi^{\alpha+1}}.
\end{split}
\end{equation*}
Thus the positive term is bounded from above by
$$ c_\alpha \omega(\xi)\bigg(\frac{2\delta^{-\alpha }}{\alpha(1-\alpha)}+  2 \xi^{-\alpha} \bigg)\omega'(\xi)\leq
\frac{c_\alpha}{\delta^{\alpha}}\frac{\omega(\xi)}{\xi^{\beta}}\frac{4(\xi+\xi^{\beta})}{\alpha(1-\alpha)}\omega'(\xi)\leq
\frac{c_\alpha \delta^{-\alpha }\gamma}{\alpha(1-\alpha)}\frac{\omega(\xi)}{\xi^{\beta}},$$
where in the first inequality we have used the fact that $\xi^{\beta-\alpha}\leq \xi+\xi^\beta$.

For the negative part, we first observe that for $\xi\geq\delta$,
\begin{equation*}
 \omega(2\xi)=\omega(\xi)+\int_{\xi}^{2\xi}\omega'(\eta)\textrm{d}
 \eta\leq\omega(\xi)+ \frac{(\log2) \gamma}{4}\leq \frac{3}{2}\omega(\xi)
\end{equation*}
under the same assumptions on $\delta$ and $\gamma$ as above. Also,
taking advantage of the concavity we obtain
$\omega(2\eta+\xi)-\omega(2\eta-\xi)\leq \omega(2\xi)$ for all
$\eta\geq\frac{\xi}{2}$. Therefore
\begin{equation*}
 c_\beta \int_{\frac{\xi}{2}}^{\infty}\frac{\omega(2\eta+\xi)-\omega(2\eta-\xi)-2\omega(\xi)}{\eta^{1+\beta}}\textrm{d}
 \eta\leq - c_\beta \frac{\omega(\xi)}{2}\int_{\frac{\xi}{2}}^{\infty}\frac{1}{\eta^{1+\beta}} \textrm{d}
 \eta \leq -\frac{ c_\beta}{2}\frac{\omega(\xi)}{\xi^{\beta}}.
\end{equation*}
But $\frac{\omega(\xi)}{\xi^{\beta}}(\frac{c_\alpha\gamma}{\delta^{\alpha}\alpha(1-\alpha)}-\frac{\nu c_\beta }{2})<0$ if
$\gamma$ is small enough \big(i.e. $\gamma<\min\{\delta, \frac{\nu\alpha(1-\alpha)  c_\beta}{ 2 c_\alpha}\delta^{\alpha} \}$\big).

Therefore, the solution $\theta(t,x)$ obeys the MOC $\omega_\lambda(\xi)$ with $\lambda$ given by \eqref{eq lambda} for all $t\in [0,T^*[$, and this directly yields
$\sup_{t\in[0,T^*[}\|\nabla\theta(t,\cdot)\|_{L^\infty}\leq \lambda$.

\section{Eventual regularity for $\alpha\in ]0,1[$ and $\beta\in ]2\alpha,\alpha+1]$}\label{sec EvtReg}
\setcounter{section}{5}\setcounter{equation}{0}

\subsection{A regularity criterion}

We first state a regularity criterion for critical and some supercritical cases concerning the H\"older continuous solutions.
\begin{proposition}\label{prop RegC}
  Let $\nu>0$, $\alpha\in]0,1[$, $\beta\in ]2\alpha, \alpha+1]$ and $\theta(t,x)$ be a smooth solution of \eqref{app-gQG}
  with $\epsilon>0$ and $\theta_0(x)\in H^m(\mathbb{R}^2)$, $m>2$. Suppose that
  $\theta(t,x)$ satisfies the bound $\norm{\theta}_{L^\infty([t_0,T]; C^{\sigma})}\leq C_0$ with $0< t_0< T<\infty$,
  $\sigma>\max\{1+\alpha-\beta, \alpha/2 \}$, and $C_0>0$ an absolute constant independent of $\epsilon$. Then for all $t\in]t_0,T]$, $s>0$ and
  $\tilde{p}>2$, we have
  $\norm{\theta(t)}_{B^s_{\tilde{p},2}}\leq f(C_0, s,\tilde{p}, t)<\infty$,
  where the function $f(C_0,s,\tilde{p}, t)$ does not depend on the value of $\epsilon$ in \eqref{app-gQG}.
\end{proposition}

\begin{proof}
We mainly follow the method from \cite{ConIW} and \cite{ConW}. Denote $\sigma_0:=\max\{\alpha+1-\beta,\frac{\alpha}{2} \}$ and we first consider $\sigma\in ]\sigma_0,1[$. From the classical $L^2$ energy method, we know the following uniform estimate
\begin{equation*}
  \norm{\theta(t)}_{L^2}^2 + 2\nu \int_0^t\|\Lambda^{\beta/2}\theta(\tau) \|_{L^2}^2\mathrm{d}\tau \leq \norm{\theta_0}_{L^2}^2, \quad t\in [0,\infty[,
\end{equation*}
Since $\sigma>\sigma_0>0$, thus by interpolation, we immediately obtain that for every $t\in [t_0,T]$ and $p\geq 2$, $\sigma_1=\sigma(1-\frac{2}{p})$
\begin{equation*}
  \|\theta(t)\|_{\dot B^{\sigma_1}_{p,\infty}} \leq \norm{\theta(t)}_{\dot B^\sigma_{\infty,\infty}}^{1-\frac{2}{p}}\norm{\theta(t)}_{L^2}^{\frac{2}{p}}\leq C_0^{1-2/p}\norm{\theta_0}_{L^2}^{4/p}.
\end{equation*}
Due to $\sigma\in]\sigma_0,1[$, we can choose
\begin{equation}
p>p_1:=\frac{2}{1-\sigma_0}
\end{equation}
such that $\sigma_1 >\sigma_0$. Then we claim that, from the a priori uniform estimate
\begin{equation}\label{eq C_1}
  \norm{\theta}_{L^\infty([t_0,T]; \dot B^{\sigma_1}_{p,\infty}\cap \dot B^{\sigma_1}_{\infty,\infty})} \leq C_1,
\end{equation}
with $C_1$ is a constant independent of $\epsilon$, we can show an improvement of the regularity that there is $\sigma_2$ chosen later such that $\sigma_2>\sigma_1$ and for every $t_1\in ]t_0,T]$
\begin{equation*}
  \sup_{t\in [t_1,T]}\|\theta(t)\|_{\dot B^{\sigma_2}_{p,\infty}\cap \dot B^{\sigma_2}_{\infty,\infty}}\leq g(t_1,\sigma_1, C_1),
\end{equation*}
where $g$ is a function given in the sequel independent of $\epsilon$.

Indeed, we apply the homogeneous dyadic operator $\dot\Delta_q$ ($q\in\mathbb{Z}$) to \eqref{app-gQG} to obtain
\begin{equation}\label{eq LocEq}
  \partial_t\dot\Delta_q\theta +\dot S_{q+1} u\cdot \nabla\dot\Delta_q\theta + \nu \Lambda^\beta\dot\Delta_q\theta -\epsilon \Delta \dot\Delta_q\theta=  \dot F_q(u,\theta),
\end{equation}
with $u=\Lambda^\alpha\mathcal{R}^{\bot}\theta$ and
\begin{equation*}
    \dot F_q(u,\theta):= \dot S_{q+1}u\cdot\nabla\dot\Delta_q \theta -\dot\Delta_q(u\cdot\nabla \theta).
\end{equation*}
Multiplying both sides of \eqref{eq LocEq} by $|\dot\Delta_q\theta|^{p-2}\dot\Delta_q \theta$, integrating over the spacial variable and using the divergence-free property and the following generalized Bernstein inequality in \cite{ChenMZ}
\begin{equation*}
  \int_{\mathbb{R}^n} \Lambda^\gamma\dot\Delta_q f(x) \; |\dot\Delta_q f|^{p-2}\dot\Delta_q f(x) \mathrm{d}x\gtrsim c 2^{q\gamma} \|\dot\Delta_q f\|_{L^p}^p,\quad  \gamma\in [0,2]
\end{equation*}
with $c$ an absolute constant independent of $q$, we have
\begin{equation*}
\begin{split}
  \frac{1}{p}\frac{d}{ dt}\| \dot\Delta_q\theta(t)\|_{L^p}^p + c\nu 2^{q\beta}\|\dot\Delta_q\theta(t)\|_{L^p}^p &\leq \Big|\int_{\mathbb{R}^2} \dot F_q(u,\theta) |\dot\Delta_q\theta|^{p-2}\dot\Delta_q\theta \mathrm{d}x\Big| \\
  & \leq \| \dot F_q(u,\theta)\|_{L^p} \| \dot \Delta_q\theta\|_{L^p}^{p-1}.
\end{split}
\end{equation*}
Thus,
\begin{equation*}
  \frac{d}{dt}\|\dot\Delta_q\theta(t) \|_{L^p} + c\nu 2^{q\beta}\|\dot \Delta_q\theta(t)\|_{L^p} \leq \| \dot F_q(u,\theta)\|_{L^p}.
\end{equation*}
It follows that for every $t\in [t_0,T]$
\begin{equation*}
  \|\dot\Delta_q\theta(t)\|_{L^p}\leq e^{-c\nu 2^{q\beta}(t-t_0)} \|\dot \Delta_q\theta(t_0)\|_{L^p} + \int_{t_0}^t e^{-c\nu2^{q\beta}(t-\tau)}\|\dot F(u,\theta)(\tau) \|_{L^p}\mathrm{d}\tau.
\end{equation*}
On the other hand, for $\dot F_q(u,\theta)$, by virtue of Bony decomposition we get
\begin{equation*}
\begin{split}
  \dot F_q(u,\theta) =&  \sum_{|k-q|\leq 1} (\dot S_{q+1}u-\dot S_{k-1}u)\cdot\nabla\dot\Delta_k\dot\Delta_q\theta -
  \sum_{|k-q|\leq 4}[\dot\Delta_q, \dot S_{k-1}u\cdot\nabla]\dot\Delta_k \theta \\
  & -\sum_{|k-q|\leq 4} \dot\Delta_q(\dot\Delta_k u\cdot\nabla \dot S_{k-1}\theta)
  - \sum_{k\geq q-3}\dot\Delta_q(\dot\Delta_k u\cdot\nabla \widetilde{\dot\Delta}_{k}\theta)
  \\ :=&  I +II+III+IV,
\end{split}
\end{equation*}
where $[A,B]:=AB-BA$ denotes the commutator operator and $\widetilde{\dot\Delta}_j:=\dot\Delta_{j-1} +\dot\Delta_j +\dot\Delta_{j+1}$.
For $I$, we directly have
\begin{equation*}
  \norm{I}_{L^p}\lesssim  2^q \|\dot\Delta_{q}\theta \|_{L^\infty} \sum_{q-2\leq q'\leq q}\|\dot\Delta_{q'}u \|_{L^p} \lesssim 2^{q(\alpha+1-2\sigma_1)}\norm{\theta}_{\dot B^{\sigma_1}_{p,\infty}} \norm{\theta}_{\dot B^{\sigma_1}_{\infty,\infty}}.
\end{equation*}
From the expression formula of $\dot\Delta_q$ and the mean value theorem, we can estimate $II$ as
\begin{equation*}
\begin{split}
  \norm{II}_{L^p} & \lesssim \sum_{|k-q|\leq 4} 2^{-q}\|\nabla\dot S_{k-1}u \|_{L^\infty} \| \nabla \dot\Delta_k \theta\|_{L^p} \\
  & \lesssim 2^{-q\sigma_1} \|\theta\|_{\dot B^{\sigma_1}_{p,\infty}}\sum_{q'\leq q+2} 2^{q'(1+\alpha-\sigma_1)} 2^{q'\sigma_1}\| \dot\Delta_{q'}\theta\|_{L^\infty} \lesssim 2^{q(1+\alpha-2\sigma_1)}\norm{\theta}_{\dot B^{\sigma_1}_{p,\infty}} \norm{\theta}_{\dot B^{\sigma_1}_{\infty,\infty}}.
\end{split}
\end{equation*}
For $III$, we obtain that for $\sigma_1<1$
\begin{equation*}
\begin{split}
  \norm{III}_{L^p}& \lesssim \sum_{|k-q|\leq 4} \|\dot\Delta_k u\|_{L^p} \|\nabla\dot S_{k-1}\theta \|_{L^\infty}
  \\ & \lesssim 2^{q(\alpha-\sigma_1)}\|\theta\|_{\dot B^{\sigma_1}_{p,\infty}}\sum_{q'\leq q+2} 2^{q'(1-\sigma_1)}2^{q'\sigma_1}
  \|\dot \Delta_{q'}\theta\|_{L^\infty}\lesssim 2^{q(1+\alpha-2\sigma_1)}\norm{\theta}_{\dot B^{\sigma_1}_{p,\infty}}
  \norm{\theta}_{\dot B^{\sigma_1}_{\infty,\infty}}.
\end{split}
\end{equation*}
For the last term, since $u$ is divergence free and $\sigma_1>\frac{\alpha}{2}$ we get
\begin{equation*}
\begin{split}
  \|IV\|_{L^p} & \lesssim 2^q\sum_{k\geq q-3} \| \dot \Delta_k u\|_{L^p}\|\widetilde{\dot\Delta}_k\theta \|_{L^\infty} \\
  & \lesssim 2^q \|\theta \|_{\dot B^{\sigma_1}_{p,\infty}}\sum_{k\geq q-4} 2^{k(\alpha-2\sigma_1)}2^{k\sigma_1}\|\dot\Delta_k \theta\|_{L^\infty} \lesssim
  2^{q(1+\alpha-2\sigma_1)}\norm{\theta}_{\dot B^{\sigma_1}_{p,\infty}} \norm{\theta}_{\dot B^{\sigma_1}_{\infty,\infty}}.
\end{split}
\end{equation*}
Hence collecting the upper estimates we infer that for every $t\in [t_0,T]$
\begin{equation*}
\begin{split}
  \|\dot\Delta_q\theta(t)\|_{L^p}& \leq e^{-c\nu 2^{q\beta}(t-t_0)} \|\dot \Delta_q\theta(t_0)\|_{L^p} +
   C_\alpha(\sigma_1)\int_{t_0}^t e^{-c\nu2^{q\beta}(t-\tau)}2^{q(1+\alpha-2\sigma_1)}\|\theta(\tau)\|_{\dot B^{\sigma_1}_{p,\infty}} \|
   \theta(\tau)\|_{\dot B^{\sigma_1}_{\infty,\infty}}\mathrm{d}\tau \\
 & \leq e^{-c\nu 2^{q\beta}(t-t_0)} \|\dot \Delta_q\theta(t_0)\|_{L^p} + \frac{C_{\alpha,\beta}(\sigma_1) (1- e^{-c\nu 2^{q\beta}(t-t_0)})}{\nu}2^{q(1+\alpha-\beta -2\sigma_1)} \|\theta\|_{L^\infty_{[t_0,t]} (\dot B^{\sigma_1}_{p,\infty}\cap \dot B^{\sigma_1}_{\infty,\infty})}^2
\end{split}
\end{equation*}
Multiplying both sides by $2^{q(2\sigma_1+\beta-\alpha-1)}$ and taking the supremum of $q$, we obtain
\begin{equation*}
 \|\theta(t)\|_{\dot B^{2\sigma_1 +\beta-\alpha-1}_{p,\infty}} \leq \sup_{q\in\mathbb{Z}} \big(e^{-c\nu 2^{q\beta}(t-t_0)} 2^{q(\sigma_1+\beta-\alpha-1)} \big)
 C_1 +  C_{\alpha,\beta}(\sigma_1) \nu^{-1} C_1^2,
\end{equation*}
where $C_1$ is from \eqref{eq C_1}.
From $\sigma_1> 1+\alpha-\beta$, we have for all $t_1\in]t_0,T]$,
\begin{equation*}
  \sup_{t\in[t_1,T]}\|\theta(t)\|_{\dot B^{2\sigma_1 +\beta-\alpha-1}_{p,\infty}} \leq \sup_{q\in\mathbb{N}} \big(e^{-c\nu 2^{q\beta}(t_1-t_0)} 2^{q(\sigma_1+\beta-\alpha-1)} \big)
   C_1 +  C_{\alpha,\beta}(\sigma_1) \nu^{-1} C_1^2:=g(t_1,\sigma_1,C_1).
\end{equation*}
Thus for all $t_1\in]t_0,T]$, we have $\theta\in L^\infty([t_1,T]; \dot B^{2\sigma_1+\beta-\alpha-1}_{p,\infty})$. Clearly $2\sigma_1+\beta-\alpha-1>\sigma_1$; in addition, from Besov embedding $ \dot B^{2\sigma_1+\beta-\alpha-1}_{p,\infty}\hookrightarrow \dot B^{2\sigma_1+\beta-\alpha-1-\frac{2}{p}}_{\infty,\infty}$, we also need
$2\sigma_1+\beta-\alpha-1-\frac{2}{p}> \sigma_1 $, which in turn leads to
\begin{equation}
  p>p_2:= \frac{2}{\sigma_1-(1+\alpha-\beta)}.
\end{equation}
Hence for $\sigma_2:=2\sigma_1+\beta-\alpha-1-\frac{2}{p}$ with $p>\max\{p_1,p_2\}$ and for all $t_1\in ]t_0,T]$, we have
\begin{equation*}
  \theta\in L^\infty([t_1,T]; \dot B^{\sigma_2}_{p,\infty}\cap \dot B^{\sigma_2}_{\infty,\infty}),
\end{equation*}
with
$$\|\theta\|_{L^\infty ([t_1,T];\dot B^{\sigma_2}_{p,\infty})}\lesssim \|\theta\|_{L^{\infty}([t_1,T];\dot B^{2\sigma_1+\beta-\alpha-1}_{p,\infty}\cap L^2)}\leq C' g(t_1,\sigma_1,C_1)$$
and
$$\|\theta\|_{L^\infty ([t_1,T];\dot B^{\sigma_2}_{\infty,\infty})}\lesssim \|\theta\|_{L^\infty([t_1,T];\dot B^{2\sigma_1+\beta-\alpha-1}_{p,\infty})}\leq C' g(t_1,\sigma_1,C_1).$$

Next we can iterate the above process through replacing $\sigma_1,\sigma_2$ by $\sigma_2,\sigma_3$ and so on;
that is, from $ \theta\in L^\infty([t_{N-1},T];\dot B^{\sigma_{N}}_{p,\infty}\cap \dot B^{\sigma_{N}}_{\infty,\infty})$, $N\geq 1$ with $t_{N-1}<t_{N}<T$, and
$$
 \|\theta\|_{L^\infty ([t_{N-1},T];\dot B^{\sigma_{N}}_{p,\infty}\cap \dot B^{\sigma_{N}}_{\infty,\infty})}\leq g_{N-1}(t_{N-1}),
$$
where $g_0(t_0)=C_1$, and
$$
g_{N-1}(t_{N-1}):=C' g(t_{N-1},\sigma_{N-1}, g_{N-2}(t_{N-2})),\quad N\geq 2,
$$
then we definitely get
$ \theta\in L^\infty([t_{N},T];\dot B^{\sigma_{N+1}}_{p,\infty}\cap \dot B^{\sigma_{N+1}}_{\infty,\infty})$ with
$$ \|\theta\|_{L^\infty([t_{N},T];\dot B^{\sigma_{N+1}}_{p,\infty}\cap \dot B^{\sigma_{N+1}}_{\infty,\infty})}\leq g_{N}(t_{N})$$
provided that $\sigma_{N}<1$ and $p>\max\{p_1,p_2\}$. Note that we do not need any additional assumption on $p$, since for every $N\geq 1$, $\sigma_{N+1}>\sigma_{N}$ and $\sigma_{N+1}+\beta -1-\alpha-\frac{2}{p}>0$. In fact, from the relation $\sigma_{N+1}=2\sigma_{N}+\beta-1-\alpha-\frac{2}{p}$, $N\geq 1$, we explicitly have
$$
\sigma_{N+1}=2^{N}(\sigma_1+\beta-1-\alpha-2/p)+ 1+\alpha+2/p-\beta,\quad N\geq 1.
$$
Due to that for some $p>\max\{p_1,p_2\}$, the fixed increment $\sigma_1+\beta-1-\alpha-\frac{2}{p}$ is positive, then $\sigma_{N+1}$ is ever increasing and it can always be attained provided that $\sigma_{N}<1$. Thus after a finite time iteration, say $N_0$ ($N_0\geq 1$), we obtain that $\sigma_{N_0+1}>1$ and
$$
 \|\theta(t)\|_{\dot B^{\sigma_{N_0+1}}_{p,\infty}\cap \dot B^{\sigma_{N_0+1}}_{\infty,\infty}}\leq g_{N_0}(t_{N_0}),\quad \forall t\in [t_{N_0},T].
$$
From $\mathcal{C}^{\sigma_{N_0+1}}= \dot B^{\sigma_{N_0+1}}_{\infty,\infty}\cap L^\infty$, we also get $\|\theta\|_{L^\infty([t_{N_0},T]; \mathcal{C}^{\sigma_{N_0+1}})}\leq C' g_{N_0}(t_{N_0})$.

We then devote to find some $p'\in ]2,\infty[$ and $\tilde{\sigma}$ such that $\tilde{\sigma}> 1+\frac{2}{p'}$ and $\theta(t)\in L^\infty ([t_{N_0},T]; B^{\tilde{\sigma}}_{p',2})$ uniformly in $\epsilon$. Let $\sigma'= \sigma_{N_0+1}(1-\frac{2}{p'})$, by interpolation we get the uniform bound that for every $t\in [t_{N_0},T]$,
\begin{equation*}
  \|\theta(t)\|_{B^{\sigma'}_{p',\infty}} \lesssim \|\theta_0\|_{L^2}+ \|\theta(t)\|_{\dot B^{\sigma_{N_0+1}}_{\infty,\infty}}^{1-\frac{2}{p'}} \|\theta_0\|_{L^2}^{\frac{2}{p'}}
  \leq C' g_{N_0}(t_{N_0}).
\end{equation*}
To make $\sigma'>1+\frac{2}{p'}$, we only need to choose some $p'$ satisfying $p'>\frac{2(\sigma_{N_0+1}+1)}{\sigma_{N_0+1}-1}$;
thus for the fixed appropriate $\sigma'$, we choose
$\tilde{\sigma}$ satisfying $1+\frac{2}{p'}<\tilde{\sigma}<\sigma'$, hence the claim follows from the continuous Besov embedding
$B^{\sigma'}_{p',\infty}\hookrightarrow B^{\tilde{\sigma}}_{p',2}$.

Now for such $\tilde{\sigma},p'$, in a similar way as obtaining \eqref{eq BmpEst1} and \eqref{eq BmpUniBd}, we get for every $t\in [t_{N_0},T]$
\begin{equation*}
\begin{split}
  \| \theta(t)\|_{B^{\tilde{\sigma}}_{p',2}}^2 + \int_{t_{N_0}}^t\|\theta(\tau)\|^2_{B^{\tilde{\sigma}+\beta/2}_{p',2}}\mathrm{d}\tau
  \leq \| \theta(t_{N_0})\|_{B^{\tilde{\sigma}}_{p',2}}^2 e^{C(t-t_{N_0})+ C \int_{t_{N_0}}^t \|\theta (\tau)\|_{\mathcal{C}^{\sigma_{N_0+1}}}^2\mathrm{d}\tau}.
\end{split}
\end{equation*}
From the equation of $\Theta_\gamma(t,x):=(t-t_{N_0})^\gamma \theta(t,x)$ ($\gamma>0$)
\begin{equation*}
  \partial_t \Theta_\gamma+ u\cdot\nabla\Theta_\gamma +\nu \Lambda^\beta \Theta_\gamma
  -\epsilon\Delta \Theta_\gamma= \gamma\Theta_{\gamma-1},\quad \Theta_\gamma(t_{N_0})=0,
\end{equation*}
with $u=\Lambda^\alpha \mathcal{R}^\bot \theta$, we infer that for every $t\in [t_{N_0},T]$
\begin{equation}\label{eq SmtEs}
  \|\Theta_\gamma(t)\|^2_{B^{\tilde{\sigma}+\gamma\beta}_{p',2}}+\int_{t_{N_0}}^t \|\Theta_\gamma(\tau)\|^2_{B^{\tilde{\sigma}+(\gamma+\frac{1}{2})\beta}_{p'2}}\mathrm{d}\tau \leq C \| \theta(t_{N_0})\|_{B^{\tilde{\sigma}}_{p',2}}^2 e^{C(\gamma+1)(t-t_{N_0} +\int_{t_{N_0}}^t \|\theta(\tau)\|_{\mathcal{C}^{\sigma_{N_0+1}}}^2\mathrm{d}\tau)}.
\end{equation}
Indeed, it reduces to consider the $\gamma=1$ case; the other cases will follow by induction or by interpolation (cf. \cite{MiaoXue}).
Similarly as \eqref{eq BmpEst1}, we have
\begin{equation*}
  \frac{d}{d t}\|\Theta_1(t)\|_{B^{\tilde{\sigma}+\beta}_{p',2}}^2 + \|\Theta_1(t)\|^2_{B^{\tilde{\sigma}+\frac{3}{2}\beta}_{p',2}}
  \leq C (1+\|\theta(t)\|_{\mathcal{C}^{\sigma_{N_0+1}}}^2) \|\Theta_1(t)\|_{B^{\tilde{\sigma}+\beta}_{p',2}}^2 + C \|\theta(t)\|_{B^{\tilde{\sigma}+\frac{1}{2}\beta}_{p',2}}^2.
\end{equation*}
Gronwall inequality ensures that
\begin{equation*}
\begin{split}
    \|\Theta_1(t)\|_{B^{\tilde{\sigma}+\beta}_{p',2}}^2 + \int_{t_{N_0}}^t\|\Theta_1(\tau)\|^2_{B^{\tilde{\sigma}+\frac{3}{2}\beta}_{p',2}}\mathrm{d}\tau
  &\leq C e^{C(t-t_{N_0})+C\int_{t_{N_0}}^t \|\theta(\tau)\|_{\mathcal{C}^{\sigma_{N_0}}}^2\mathrm{d}\tau}\int_{t_{N_0}}^t \|\theta(\tau)\|_{B^{\tilde{\sigma}+\frac{\beta}{2}}_{p',2}}^2\mathrm{d}\tau \\
 & \leq C \|\theta(t_{N_0})\|_{B^{\tilde{\sigma}}_{p',2}}^2 e^{2C(t-t_{N_0})+2C\int_{t_{N_0}}^t \|\theta(\tau)\|_{\mathcal{C}^{\sigma_{N_0+1}}}^2\mathrm{d}\tau}.
\end{split}
\end{equation*}
\eqref{eq SmtEs} means that for every $s>0$ and $t\in ]t_{N_0},T]$, we have the uniform bound of the norm $\|\theta(t)\|_{B^{s}_{p',2}}$.
Therefore the conclusion follows from Besov embedding or interpolation (with $L^2$) and the fact that
$t_1<\cdots<t_{N_0}$ are chosen arbitrarily in $]t_0,T]$.

\end{proof}

\subsection{Proof of eventual regularity for $\alpha\in]0,1[$ and $\beta \in ]2\alpha,\alpha+1]$}

Indeed, according to Proposition \ref{prop RegC} and a standard process of approximation, it suffices to show the following uniform H\"older estimate of the solution.
\begin{theorem}\label{thm HolEs}
Let $\nu>0$, $\alpha\in ]0,1[$, $\beta\in ]2\alpha, \alpha+1]$ and $\theta(t,x)$ be the global smooth solution of \eqref{app-gQG} with $\epsilon>0$ and $\theta_0(x)\in H^m(\mathbb{R}^2)$, $m>2$. Then for every $\gamma\in ]\max\{\alpha+1-\beta,\alpha/2 \}, 1[$, there exists a time $T$ depending only on $ \alpha,\beta,\gamma,\nu, \|\theta_0\|_{L^\infty}$ such that $ \|\theta(t,\cdot)\|_{\mathcal{C}^\gamma}$ is uniformly bounded with respect to $\epsilon$ for all $t\geq T$.
\end{theorem}

We shall use the nonlocal maximum principle method to prove Theorem \ref{thm HolEs} (cf.\cite{Kis}). We first claim that some stationary moduli of continuity which imply the uniform H\"older estimates can be preserved by the evolution of the critical and some supercritical equation \eqref{app-gQG}.
\begin{proposition}\label{prop sMOC}
Let $\nu>0$, $\alpha\in ]0,1[$, $\beta\in ]2\alpha, \alpha+1]$, $\gamma\in ]\max\{\alpha+1-\beta,\alpha/2 \}, 1[$ and $\theta(t,x)$ be the global smooth solution of \eqref{app-gQG} with $\epsilon>0$ and $\theta_0(x)\in H^m(\mathbb{R}^2)$, $m>2$. For every $H,\delta\in\mathbb{R}^+$, set
\begin{equation}\label{sMOC}
  \omega(\xi)=
  \begin{cases}
    (H/\delta^\gamma)\xi^\gamma, &\quad \textrm{if} \quad \xi\in ]0,\delta],\\
    H,&\quad \textrm{if}\quad \xi\in ]\delta,\infty[.
  \end{cases}
\end{equation}
Suppose that $H\geq 2(1+c_1)\|\theta_0\|_{L^\infty}$ is satisfied for some $c_1\in ]0,1/2]$ and the initial data $\theta_0(x)$ obeys $\omega(\xi)$. Then there exists an absolute constant $C_1=C_1(\alpha,\beta,\gamma,\nu)$ such that if $H \leq C_1
\delta^{\alpha+1-\beta}$, the solution $\theta(t,x)$ obeys $\omega(\xi)$ for all $t > 0$, independently of $\epsilon$.
\end{proposition}

The proof of Proposition \ref{prop sMOC} is placed in the sequel of this subsection. Clearly, the function $\omega(\xi)$ from \eqref{sMOC}
is a modulus of continuity satisfying $\omega(0+)= 0$ and $\omega'(0+)=+\infty$ (condition (b)). Moreover,
if $\theta(x,t)$ obeys some $\omega(\xi)$ for a time period $I$, it implies the uniform H\"older estimate
$\|\theta(t,\cdot)\|_{\mathcal{C}^\gamma}\leq \frac{H}{\delta^\gamma}+ \|\theta_0\|_{L^\infty}$ for all $t\in I$.
Hence for Theorem \ref{thm HolEs}, it remains to show that after a finite time $T$ independent of $\epsilon$,
$\theta(t,x)$ obeys some moduli from \eqref{sMOC} for all $t\geq T$.

Here we note that, due to the restriction $2(1+c_1)\|\theta_0\|_{L^\infty}\leq H\leq C_1 \delta^{\alpha+1-\beta}$, not every initial data will obey some such moduli. Observe that for $|x-y|\leq \delta$, we have
\begin{equation*}
  |\theta_0(x)-\theta_0(y)|\leq \|\theta_0\|_{\mathcal{C}^\gamma}|x-y|^\gamma = \frac{\|\theta_0\|_{\mathcal{C}^\gamma}\delta^\gamma}{H}\omega(|x-y|)\leq
  \frac{\|\theta_0\|_{\mathcal{C}^\gamma}\delta^\gamma}{2(1+c_1)\|\theta_0\|_{L^\infty}} \omega(|x-y|);
\end{equation*}
while for $|x-y|>\delta$, we get $|\theta_0(x)-\theta_0(y)|\leq 2\|\theta_0\|_{L^\infty}<H$. Thus $\delta\in\mathbb{R}^+$ should satisfy $C_1 \delta^{\alpha+1-\beta}> 2(1+c_1)\|\theta_0\|_{L^\infty}$ and $\delta^\gamma<\frac{2(1+c_1)\|\theta_0\|_{L^\infty}}{\|\theta_0\|_{\mathcal{C}^\gamma}}$. Clearly, we can not choose some appropriate
$\delta$ so that the two conditions simultaneously hold for all the initial data. Even we use the decay estimate of $L^\infty$ norm
(Lemma \ref{lem LinfDC}) and consider the solution $\tilde{\theta}(t,x)=\theta(t+T,x)$ instead of $\theta(t,x)$, we still can not find some suitable $\delta$
uniform in $\epsilon$. However, corresponding to the eventual nature of Theorem \ref{thm HolEs}, we shall show that the solution will obey some suitable moduli from \eqref{sMOC} after a finite time $T$ independent of $\epsilon$, so that Proposition \ref{prop sMOC} can be applied to the time-translated solution $\tilde{\theta}(t,x)=\theta(t+T,x)$.

Before stating the next key Lemma, we introduce a variant family of moduli
\begin{equation}\label{evtMOC}
  \omega(\xi,\xi_0)=\begin{cases}
    \gamma H \delta^{-\gamma}\xi_0^{\gamma-1}\xi +(1-\gamma) H \delta^{-\gamma}\xi_0^\gamma,&\quad \textrm{if} \quad \xi \in ]0,\xi_0], \\
    (H/\delta^\gamma)\xi^\gamma,  &\quad \textrm{if}\quad \xi\in ]\xi_0, \delta], \\
    H,  &\quad \textrm{if}\quad \xi\in ]\delta,\infty[,
  \end{cases}
\end{equation}
with $0< \xi_0\leq\delta$, $\gamma\in ]0,1[ $ and $H,\delta\in ]0,\infty[$. Compared with $\omega(\xi)$ from \eqref{sMOC}, the only difference lies in that when $\xi\in]0,\xi_0] $, $\omega(\xi)$ is replaced by the tangent line of $\omega(\xi)$ at $\xi=\xi_0$. Clearly $\omega(\xi,\xi_0)$ is a modulus of continuity satisfying that $\omega(0+,\xi_0)>0$ for all $\xi_0>0$ (condition (a)) and $\omega(\xi,0+)=\omega(\xi)$. Also, any bounded initial data $\theta_0(x)$ obeys $\omega(\xi,\delta)$ if $2\|\theta_0\|_{L^\infty}\leq \omega(0+,\delta)= H(1-\gamma)$ is satisfied. Then we claim that

\begin{lemma}\label{lem evtMoc}
Let $\nu>0$, $\alpha\in ]0,1[$, $\beta\in ]2\alpha, \alpha+1]$, $\gamma\in ]\max\{\alpha+1-\beta,\alpha/2 \}, 1[$
and $\theta(t,x)$ be the global smooth solution of \eqref{app-gQG} with $\epsilon>0$ and $\theta_0(x)\in H^m(\mathbb{R}^2)$, $m>2$.
Suppose that $H\geq 2(1+c_1)\|\theta_0\|_{L^\infty}$ is satisfied for some $c_1\in ]0,1/2]$ and the initial data $\theta_0(x)$ obeys $\omega(\xi,\delta)$. Then there exist two positive absolute constants $C_1=C_1(\alpha,\beta,\gamma,\nu)$ and $C_2=C_2(\beta,\gamma,\nu)$ such that if $ H\leq C_1\delta^{\alpha+1-\beta}$ and $\xi_0(t)$ is the solution of
\begin{equation*}
  \frac{d}{dt}\xi_0(t)= -C_2 \xi_0(t)^{1-\beta},\quad \xi_0(0)=\delta,
\end{equation*}
then the solution $\theta(t,x)$ obeys the modulus of continuity $\omega(\xi,\xi_0(t))$ for all $t$ satisfying $\xi_0(t)> 0$, independent of $\epsilon$.
\end{lemma}

The proof of Lemma \ref{lem evtMoc} is also placed in the sequel of this subsection. Note that at the current state $C_1$ and $C_2$
are two fixed absolute constants and the two $C_1$ in Lemma \ref{lem evtMoc} and Proposition \ref{prop sMOC} can be chosen to be identical.

Since $\xi_0(t)=(\delta^\beta-C_2\beta t)^{1/\beta}$, we know that after a finite time $T_0=\delta^\beta/(C_2\beta)$, $\theta(T_0-,x)=\theta(T_0,x)$ will obey the MOC $\omega(\xi,0+)=\omega(\xi)$ and thus Proposition \ref{prop sMOC} can be applied to the solution $\tilde{\theta}(t,x):=\theta(t+T_0,x)$. Thus to prove Theorem \ref{thm HolEs}, it suffices to find some suitable $H,\delta,c_1$ under the conditions that
$$2\|\theta_0\|_{L^\infty}\leq H(1-\gamma)\quad \textrm{and} \quad 2(1+c_1)\|\theta_0\|_{L^\infty}\leq H\leq C_1\delta^{\alpha+1-\beta}.$$
The supercritical case is direct, since $\alpha+1>\beta$, we only need to choose $c_1=\min\{\frac{1}{1-\gamma}-1,\frac{1}{2}\}$ and $\delta$ large enough such that $\delta^{\alpha+1-\beta}>\frac{ 2\|\theta_0\|_{L^\infty}}{C_1(1-\gamma)}$, then some appropriate $H$ satisfying $\frac{2}{1-\gamma}\|\theta_0\|_{L^\infty}\leq H \leq C_1 \delta^{\alpha+1-\beta}$ can also be chosen.
For the critical case, we shall use the decay estimate Lemma \ref{lem LinfDC} to obtain that after a finite time independent of $\epsilon$, say $\tilde{T}$, $\|\theta(\tilde{T},x)\|_{L^\infty}< C_1(1-\gamma)/2$, then we only need to choose the same $c_1$ as in the supercritical case, some $H$ satisfying $\frac{2}{1-\gamma}\|\theta(\tilde{T})\|_{L^\infty}\leq H\leq C_1$ and replace $\theta(t,x)$ by $\bar{\theta}(t,x)=\theta(t+\tilde{T},x)$.

Now, we are devoted to the proof of Proposition \ref{prop sMOC} and Lemma \ref{lem evtMoc}. We shall apply the general criterion Proposition \ref{prop GC}, and due to that the solution $\theta(t,x)$ (or the translated solution $\tilde{\theta}(t,x)$)
and the moduli $\omega(\xi)$, $\omega(\xi,\xi_0(t))$ satisfy the conditions needed,
it reduces to check the inequality \eqref{eq keyGC}.
Thus first we should know the expression of $\Omega(t,\xi)$ in \eqref{eq keyGC}. In fact, from Lemma 3.2 in \cite{MiaoXue}, we have had a rough
estimate: let $\theta(t,x)$ obeys $\omega(t,\xi)$, then $u= \Lambda^\alpha \mathcal{R}^\bot \theta$ satisfies that for every $x,y\in\mathbb{R}^2$
\begin{equation}\label{eq uEsRo}
  |u(t,x)-u(t,y)|\leq c_\alpha \Big( \int_0^\xi \frac{\omega(t,\eta)}{\eta^{1+\alpha}}\mathrm{d}\eta + \xi\int_\xi^\infty \frac{\omega(t,
\eta)}{\eta^{2+\alpha}}\mathrm{d}\eta\Big)
\end{equation}
with $\xi=|x-y|$ and $c_\alpha$ an absolute constant. The RHS of the expression has very strong singularity at the origin, and especially for the moduli $\omega(\xi,\xi_0(t))$ with $\omega(0+,\xi_0(t))>0$, it clearly diverges. However, similarly as treating the corresponding case of quasi-geostrophic equation in \cite{Kis},
we can rely on an improved estimate of the dissipative term in the breakdown scenario to overcome the difficulty. Indeed, we have
\begin{lemma}\label{lem DispImEs}
Let $\nu>0$, $\alpha\in ]0,1[$, $\beta\in ]2\alpha, \alpha+1]$ and $\theta(t,x)$ be the global smooth solution of \eqref{app-gQG} with $\epsilon>0$ and
$\theta_0(x)\in H^m(\mathbb{R}^2)$, $m>2$. Assume that the modulus of continuity $\omega(t,\xi)$ satisfies the assumptions in Proposition \ref{prop GC} and $\theta_0(x)$ obeys $\omega(0,\xi)$. Also suppose that $T_*$ is the minimal time that the MOC $\omega$ is lost, and $x,y\in\mathbb{R}^2$, $x\neq y$ are two points
satisfying
\begin{equation*}
  \theta(T_*,x)-\theta(T_*,y)= \omega(T_*,\xi), \quad \textrm{with}\;\; \xi:=|x-y|.
\end{equation*}
Then the following statements are true.
\begin{enumerate}[(1)]
  \item For $x_0:=(\frac{\xi}{2},0)$, $y_0:=(-\frac{\xi}{2},0)$, there exist a unique rotating transform $\rho$ and a unique vector $a\in\mathbb{R}^2$ such that \begin{equation*}
        x=\rho x_0-a,\quad y=\rho y_0-a.
      \end{equation*}
  \item We have
  \begin{equation}\label{eq DispImEs}
  -\Lambda^\beta\theta(T_*,x)+\Lambda^\beta\theta(T_*,y)\leq \Upsilon_\beta(T_*,\xi)+ \Upsilon_\beta^{\bot}(T_*,\xi).
\end{equation}
Here $\Upsilon_\beta^\bot$ is bounded by
\begin{equation}\label{eq UpsBot}
  \Upsilon_\beta^\bot(T_*,\xi)\leq -C\int_{B^+_{r_0\xi}(x_0) }\frac{f_{\rho,a}(T_*,\eta,\mu)}{|x_0-(\eta,\mu)|^{2+\beta}}\mathrm{d}\eta\mathrm{d}\mu,
\end{equation}
where
\begin{equation}\label{eq frhoa}
  f_{\rho,a}(t,\eta,\mu) := 2\omega(t,2\eta)- \theta_{\rho,a}(t,\eta,\mu)+\theta_{\rho,a}(t,-\eta,\mu)-\theta_{\rho,a}(t,\eta,-\mu)+\theta_{\rho,a}(t,-\eta,-\mu),
\end{equation}
and
\begin{equation}\label{eq thtarhoa}
  \theta_{\rho,a}(t,\eta,\mu):= \theta(t, \rho\cdot(\eta,\mu)-a),
\end{equation}
and
\begin{equation}
  B^+_{r_0\xi}(x_0):=\{(\eta,\mu)\in\mathbb{R}^2 : |x_0-(\eta,\mu)|\leq r_0\xi,\; \mu>0\}.
\end{equation}
In the above $C,r_0>0$ are absolute constants that may depend on $\beta$.
\end{enumerate}

\end{lemma}

We place the proof of this Lemma in the next subsection.

Note that $\Upsilon_\beta^\bot(T_*,\xi)\leq 0$, then according to the structure of the nonlinear term, we can control
the strong singularity near the origin in terms of $-\Upsilon_\beta^\bot$. Precisely, we have
\begin{lemma}\label{lem nLinImEs}
Under the assumptions of Lemma \ref{lem DispImEs}, and let $u=\Lambda^\alpha\mathcal{R}^\bot\theta$, $\ell:=\frac{x-y}{|x-y|}$. Then we have
\begin{equation}\label{eq nLinImEs1}
  |(u(T_*,x)-u(T_*,y))\cdot\ell|\leq \Omega(T_*,\xi),
\end{equation}
with
\begin{equation}\label{eq OmeIm}
  \Omega(t,\xi)= A \Big( -\xi^{\beta-\alpha}\Upsilon^\bot_\beta(t,\xi)+ \xi\int_\xi^\infty \frac{\omega(t,\xi)}{\eta^{2+\alpha}}\mathrm{d}\eta +\xi^{-\alpha}\omega(t,\xi)\Big),
\end{equation}
where $A $ is an absolute constant that may depend on $\alpha,\beta$.
In particular, if $\omega(t,\xi)$ is $\omega(\xi)$ from \eqref{sMOC} or $\omega(\xi,\xi_0(t))$ from \eqref{evtMOC}, we also get
\begin{equation}\label{eq OmeIm1}
  \Omega(t,\xi)= A \Big( -\xi^{\beta-\alpha}\Upsilon^\bot_\beta(t,\xi) +\xi^{-\alpha}\omega(t,\xi)\Big).
\end{equation}
\end{lemma}
We also place the proof of the Lemma in the next subsection.

Now based on the above two Lemmas, we prove Proposition \ref{prop sMOC} and Lemma \ref{lem evtMoc} as follows.
\begin{proof}[Proof of Proposition \ref{prop sMOC}]
It suffices to check the following inequality
\begin{equation}\label{eq keyGCIm1}
  \Omega(\xi)\omega'(\xi)+\nu \Upsilon_\beta(\xi)+\nu\Upsilon_\beta^\bot(\xi)<0,\quad \forall \xi>0,
\end{equation}
where $\omega,\Omega,\Upsilon_\beta,\Upsilon_\beta^\bot$ are defined by \eqref{sMOC}, \eqref{eq OmeIm1}, \eqref{eq GC-Ups} and \eqref{eq UpsBot} respectively and we have omitted the time variable due to the time independence.
The case $\xi>\delta$ is obvious from $\omega(\xi)=H$. While when $\xi\in ]0,\delta]$, we only need to prove
\begin{equation*}
  A \xi^{-\alpha}\omega(\xi)\omega'(\xi)+\nu\Upsilon_\beta(\xi)<0, \quad\textrm{and}\quad -A\xi^{\beta-\alpha}\omega'(\xi)\Upsilon_\beta^\bot(\xi)+ \nu \Upsilon_\beta^\bot(\xi)\leq 0.
\end{equation*}
From $\omega(\xi)= \frac{H}{\delta^\gamma} \xi^\gamma$, $\omega'(\xi)=\gamma\frac{H}{\delta^\gamma}\xi^{\gamma-1}$, $\omega''(\xi)=-\gamma(1-\gamma)\frac{H}{\delta^\gamma}\xi^{\gamma-2}$, and
\begin{equation*}
  \Upsilon_\beta(\xi)\leq c_\beta\int_0^{\frac{\xi}{2}}\frac{\omega(\xi+2\eta)+\omega(\xi-2\eta)-2\omega(\xi)}{\eta^{1+\beta}}\mathrm{d}\eta\leq c_\beta'\xi^{2-\beta}\omega''(\xi),
\end{equation*}
we have to call for that
\begin{equation*}
  A \gamma \frac{H^2}{\delta^{2\gamma}}\xi^{2\gamma-\alpha-1} - \nu c_\beta' \gamma(1-\gamma)\frac{H}{\delta^\gamma}\xi^{\gamma-\beta}=
  A \gamma \frac{H}{\delta^\gamma} \xi^{\gamma-\beta}\Big(\frac{H}{\delta^{\alpha+1-\beta}}\big(\frac{\xi}{\delta}\big)^{\gamma+\beta-1-\alpha} -\frac{\nu c_\beta' (1-\gamma)}{A}\Big)<0,
\end{equation*}
and
\begin{equation*}
  -A \gamma\frac{H}{\delta^\gamma}\xi^{\beta-\alpha+\gamma-1}\Upsilon^\bot_\beta(\xi) +\nu \Upsilon^\bot_\beta(\xi)=\Big(-A \gamma\frac{H}{\delta^{\alpha+1-\beta}} \big(\frac{\xi}{\delta} \big)^{\gamma+\beta-\alpha-1}+ \nu \Big)\Upsilon_\beta^\bot(\xi)\leq 0.
\end{equation*}
Due to $\xi\in]0,\delta]$, $\gamma>\alpha+1-\beta$ and $\Upsilon_\beta^\bot\leq 0$, both inequalities can be guaranteed by the condition $\frac{H}{\delta^{\alpha+1-\beta}}\leq C_1$ with $C_1$  some positive number satisfying $C_1< \min\{\nu c'_\beta(1-\gamma)/A,\nu/(A\gamma)\}$.
\end{proof}

\begin{proof}[Proof of Lemma \ref{lem evtMoc}]
It suffices to check the following inequality that for all $\xi>0$ and $t>0$
\begin{equation}\label{eq keyGCIm2}
  -\partial_t\omega(\xi,\xi_0(t))+ \Omega(t,\xi)\partial_\xi\omega(\xi,\xi_0(t)) + \nu\Upsilon_\beta(t,\xi)+ \nu\Upsilon_\beta^\bot(t,\xi)<0,
\end{equation}
where $\omega(t,\xi)=\omega(\xi,\xi_0(t))$ from \eqref{evtMOC} and
$\Omega,\Upsilon_\beta,\Upsilon_\beta^\bot$ are given by \eqref{eq OmeIm1}, \eqref{eq GC-Ups} and \eqref{eq UpsBot}
respectively. The case $\xi>\delta$ is immediate due to $\omega(\xi,\xi_0(t))=H$. For the case $\xi\in ]\xi_0(t),\delta]$, since $\omega(\xi,\xi_0(t))= \frac{H}{\delta^\gamma}\xi^\gamma$, almost parallel to the proof of Proposition \ref{prop sMOC}, we know that \eqref{eq keyGCIm2} holds for sufficiently small $C_1$.
We then consider the case $\xi\in ]0,\xi_0(t)]$, and it only needs to verify the
following estimates
\begin{equation}\label{eq evtMocPf1}
  -\partial_t \omega(\xi,\xi_0(t)) +(\nu/2)\Upsilon_\beta(t,\xi)<0,
\end{equation}
and
\begin{equation}\label{eq evtMocPf2}
  A\xi^{-\alpha}\omega(\xi,\xi_0(t))\partial_\xi\omega(\xi,\xi_0(t))+(\nu/2) \Upsilon_\beta(t,\xi)<0,
\end{equation}
and
\begin{equation}\label{eq evtMocPf3}
  -A\xi^{\beta-\alpha}\partial_{\xi}\omega(\xi,\xi_0(t))\Upsilon_\beta^\bot(t,\xi) + \nu\Upsilon^\bot_\beta(t,\xi)\leq 0.
\end{equation}
From $\omega(\xi,\xi_0(t))=\gamma \frac{H}{\delta^\gamma}\xi_0^{\gamma-1}\xi + (1-\gamma)\frac{H}{\delta^\gamma}\xi_0^\gamma$ for $\xi\in ]0,\xi_0(t)]$ and $\xi_0'(t)\leq 0$, we get
\begin{equation*}
  -\partial_t\omega(\xi,\xi_0(t))= \gamma(1-\gamma)\frac{H}{\delta^\gamma}\big(\xi_0^{\gamma-2}\xi -\xi_0^{\gamma-1} \big)\xi_0'\leq
  - \gamma(1-\gamma)\frac{H}{\delta^\gamma}\xi_0^{\gamma-1}\xi_0',
\end{equation*}
and
\begin{equation*}
  \omega(\xi,\xi_0(t))\partial_\xi\omega(\xi,\xi_0(t))= \gamma\frac{H^2}{\delta^{2\gamma}} \xi_0^{2\gamma-2}(\gamma\xi+(1-\gamma)\xi_0)\leq
  \gamma\frac{H^2}{\delta^{2\gamma}} \xi_0^{2\gamma-1}.
\end{equation*}
For the dissipation term $\Upsilon_\beta$, we obtain
\begin{equation*}
  \Upsilon_\beta(t,\xi)\leq c_\beta \int_{\frac{\xi}{2}}^\infty \frac{\omega(\xi+2\eta,\xi_0)-\omega(\xi-2\eta,\xi_0)-2\omega(\xi,\xi_0)}{\eta^{1+\beta}}
  \mathrm{d}\eta \leq -c_\beta' \frac{\omega(0,\xi_0)}{\xi^{\beta}},
\end{equation*}
where we have used the inequality $\omega(\xi+2\eta,\xi_0)-\omega(\xi-2\eta,\xi_0)-2\omega(\xi,\xi_0)\leq -2\omega(0,\xi_0)$. Then by virtue of $\xi_0'=-C_2\xi_0^{1-\beta}$
and $\omega(0,\xi_0)=(1-\gamma)\frac{H}{\delta^\gamma}\xi_0^\gamma$, the LHS of \eqref{eq evtMocPf1} can be bounded by
\begin{equation*}
  - \gamma(1-\gamma)\frac{H}{\delta^\gamma}\xi_0^{\gamma-1}\xi_0'-\frac{\nu}{2} c_\beta'
  (1-\gamma)\frac{H}{\delta^\gamma}\frac{\xi_0^{\gamma}}{\xi^\beta}=(1-\gamma)\frac{H}{\delta^\gamma}\xi_0^{\gamma-\beta}\Big(\gamma C_2-\frac{\nu c_\beta'}{2}
  \big( \frac{\xi_0}{\xi}\big)^\beta  \Big).
\end{equation*}
Clearly, the last expression is strictly negative if $C_2$ is sufficiently small (i.e. $C_2< \frac{ c_\beta'\nu }{2\gamma}$). Also, the LHS of \eqref{eq evtMocPf2} and
\eqref{eq evtMocPf3} can be respectively bounded by
\begin{equation*}
  A\gamma\frac{H^2}{\delta^{2\gamma}} \xi_0^{2\gamma-1}\xi^{-\alpha}-\frac{\nu}{2}c'_\beta (1-\gamma)\frac{H}{\delta^\gamma}\frac{\xi_0^{\gamma}}{\xi^\beta}=
  \frac{H}{\delta^\gamma}\frac{\xi_0^{\gamma}}{\xi^\beta}\Big(A\gamma \frac{H}{\delta^{\alpha+1-\beta}}\big(\frac{\xi_0}{\delta}\big)^{\gamma+\beta-\alpha-1}\big(\frac{\xi}{\xi_0}\big)^{\beta-\alpha} -\frac{\nu}{2}c'_\beta (1-\gamma)\Big),
\end{equation*}
and
\begin{equation*}
  -A\gamma \frac{H}{\delta^\gamma} \xi_0^{\gamma-1}\xi^{\beta-\alpha}\Upsilon_\beta^\bot(t,\xi)+\nu\Upsilon_\beta^\bot(t,\xi)=\Big(-A\gamma
  \frac{H}{\delta^{\alpha+1-\beta}}\big(\frac{\xi_0}{\delta}\big)^{\gamma+\beta-\alpha-1}\big(\frac{\xi}{\xi_0}\big)^{\beta-\alpha} +\nu \Big)\Upsilon_\beta^\bot(t,\xi).
\end{equation*}
Due to $\gamma>\alpha+1-\beta$ and $\beta>\alpha$, the last expressions of both formulae are (strictly) negative if $C_1$ is sufficiently small (i.e. $C_1<\min\{\frac{\nu c_\beta'(1-\gamma)}{2A\gamma},\frac{\nu}{A\gamma}\}$).
\end{proof}

\subsection{Proof of two technical Lemmas}
In this subsection we are dedicated to the proof Lemma \ref{lem DispImEs} and Lemma \ref{lem nLinImEs}.

\begin{proof}[Proof of Lemma \ref{lem DispImEs}]
Statement (1) is obvious. We then consider (2). First denote by $\mathcal{P}_{h,n}^\beta(x)$ the $n$-dimensional kernel of the operator $e^{-h\Lambda^\beta}$, indeed from \cite{CafS} we have the following explicit expression that for every $\beta\in ]0,2[$
\begin{equation}\label{eq PhnBeta}
  \mathcal{P}_{h,n}^\beta(x)= c(\beta,n)\frac{h}{(|x|^2+\beta^2h^{\frac{2}{\beta}})^{\frac{n+\beta}{2}}}.
\end{equation}
Thus we have
\begin{equation*}
  -\Lambda^\beta\theta(T_*,x)+\Lambda^\beta\theta(T_*,y)= \lim_{h\rightarrow 0}\frac{1}{h}\Big(\mathcal{P}_{h,2}^\beta*\theta(T_*,x) -\mathcal{P}_{h,2}^\beta*\theta(T_*,y) -\theta(T_*,x)+\theta(T_*,y) \Big).
\end{equation*}
Since $\theta(T_*,x)-\theta(T_*,y)=\omega(T_*,\xi)$, it only needs to estimate the difference of the remaining terms. Also for brevity we omit the time variable $T_*$ in the sequel. Then
\begin{eqnarray}
  \mathcal{P}_{h,2}^\beta*\theta(x) -\mathcal{P}_{h,2}^\beta*\theta(y) = \int_{\mathbb{R}^2}\big(\mathcal{P}_{h,2}^\beta(\rho x_0-a -z)
  -\mathcal{P}_{h,2}^\beta(\rho y_0-a-z)\big)\theta(z)\mathrm{d}z &\nonumber\\
   = \int_{\mathbb{R}^2} \big(\mathcal{P}_{h,2}^\beta(x_0-z)-\mathcal{P}_{h,2}^\beta(y_0-z)\big)\theta_{\rho,a}(z)\mathrm{d}z& \nonumber\\
  = \iint_{\mathbb{R}^2} \big(\mathcal{P}_{h,2}^\beta(\frac{\xi}{2}-\eta,\mu)-\mathcal{P}_{h,2}^\beta(-\frac{\xi}{2}-\eta,\mu)
  \big)\theta_{\rho,a}(\eta,\mu)\mathrm{d}\eta\mathrm{d}\mu & \nonumber\\
   =\int_{\mathbb{R}}\mathrm{d}\mu \int_0^\infty\big(\mathcal{P}_{h,2}^\beta(\frac{\xi}{2}-\eta,\mu)-\mathcal{P}_{h,2}^\beta(\frac{\xi}{2}+\eta,\mu)
  \big)(\theta_{\rho,a}(\eta,\mu)-\theta_{\rho,a}(-\eta,\mu))\mathrm{d}\eta & \nonumber\\
  = \int_{\mathbb{R}}\mathrm{d}\mu \int_0^\infty\big(\mathcal{P}_{h,2}^\beta(\frac{\xi}{2}-\eta,\mu)-\mathcal{P}_{h,2}^\beta(\frac{\xi}{2}+\eta,\mu)
  \big)\omega(2\eta)\mathrm{d}\eta + & \nonumber\\
  \int_{\mathbb{R}}\mathrm{d}\mu \int_0^\infty\big(\mathcal{P}_{h,2}^\beta(\frac{\xi}{2}-\eta,\mu)-\mathcal{P}_{h,2}^\beta(\frac{\xi}{2}+\eta,\mu)
  \big)(\theta_{\rho,a}(\eta,\mu)-\theta_{\rho,a}(-\eta,\mu)-\omega(2\eta))\mathrm{d}\eta & \label{eq UpsBetBot} \\
  :=\Upsilon_{\beta,h}(\xi) + \Upsilon_{\beta,h}^\bot(\xi).&
\end{eqnarray}
In a similar way as treating the corresponding part in \cite{MiaoXue} or \cite{Kis}, we know
\begin{equation*}
  \lim_{h\rightarrow 0}\frac{1}{h}(\Upsilon_{\beta,h}(\xi) -\omega(\xi))= \Upsilon_{\beta}(\xi),
\end{equation*}
with $\Upsilon_\beta(\xi)$ is given by \eqref{eq GC-Ups}. Now we estimate $\lim_{h\rightarrow 0}\frac{1}{h}\Upsilon_{\beta,h}^\bot(\xi)$.
Clearly the integrand in \eqref{eq UpsBetBot} is always negative, and for the integral we only consider a small part that is near the dangerous point $x_0=(\frac{\xi}{2},0)$.
In fact, from \eqref{eq PhnBeta} and the fact that $h$ is arbitrarily small, there exists a small universal number $r_0\in ]0,\frac{1}{4}[$ such that for every $z\in B_{r_0\xi}(x_0)\subset \mathbb{R}^2$, we have
\begin{equation*}
  \mathcal{P}_{h,2}^\beta(x_0-z)-\mathcal{P}_{h,2}^{\beta}(y_0-z)\geq \frac{1}{2}\mathcal{P}_{h,2}^\beta(x_0-z).
\end{equation*}
Thus we get
\begin{align}
  &\Upsilon_\beta^\bot(\xi) =\lim_{h\rightarrow 0}\frac{1}{h}\Upsilon_{\beta,h}^\bot(\xi) \nonumber\\
  =&  \lim_{h\rightarrow 0}\frac{1}{h} \int_{\mathbb{R}}\mathrm{d}\mu
   \int_0^\infty\big(\mathcal{P}_{h,2}^\beta(\frac{\xi}{2}-\eta,\mu)-\mathcal{P}_{h,2}^\beta(\frac{\xi}{2}+\eta,\mu)
  \big)(\theta_{\rho,a}(\eta,\mu)-\theta_{\rho,a}(-\eta,\mu)-\omega(2\eta))\mathrm{d}\eta \nonumber\\
   \leq& \lim_{h\rightarrow 0}\frac{1}{2h} \iint_{B_{r_0\xi}(x_0)}
  \mathcal{P}_{h,2}^\beta(\frac{\xi}{2}-\eta,\mu)(\theta_{\rho,a}(\eta,\mu)-\theta_{\rho,a}(-\eta,\mu)-\omega(2\eta))\mathrm{d}\eta \mathrm{d}\mu \nonumber\\
  \leq&  c(\beta)\iint_{B_{r_0\xi}(x_0)}
  \frac{\theta_{\rho,a}(\eta,\mu)-\theta_{\rho,a}(-\eta,\mu)-\omega(2\eta)}{|x_0-(\eta,\mu)|^{2+\beta}}\mathrm{d}\eta\mathrm{d}\mu \nonumber \\
  =&  - c(\beta)\iint_{B^+_{r_0\xi}(x_0)}
  \frac{2\omega(2\eta)-\theta_{\rho,a}(\eta,\mu)+ \theta_{\rho,a}(-\eta,\mu)-\theta_{\rho,a}(\eta,-\mu)+
  \theta_{\rho,a}(-\eta,-\mu)}{|x_0-(\eta,\mu)|^{2+\beta}}\mathrm{d}\eta\mathrm{d}\mu. \label{eq UpsBetBot1}
\end{align}
We note that although the denominator of \eqref{eq UpsBetBot1}  contains the non-integrable singularity, the whole integral is still absolutely integrable due to the cancelation in the
numerator. More precisely, from $\theta_{\rho,a}(x_0)-\theta_{\rho,a}(y_0)=\omega(\xi)$, we know
$\partial_1\theta_{\rho,a}(x_0)=\partial_1\theta_{\rho,a}(y_0)=\omega'(\xi)$ and $\partial_2\theta_{\rho,a}(x_0)=\partial_2\theta_{\rho,a}(y_0)=0$. Thus
according to Taylor formula, we further obtain that for $(\eta,\mu)\in B_{r_0\xi}(x_0)$
\begin{equation*}
  \theta_{\rho,a}(\eta,\mu)-\theta_{\rho,a}(\frac{\xi}{2},0)\geq \partial_1\theta_{\rho,a}(\frac{\xi}{2},0)(\eta-\frac{\xi}{2}) -\frac{\|\nabla^2\theta_{\rho,a}\|_{L^\infty}}{2}
  ((\eta-\frac{\xi}{2})^2 + \mu^2),
\end{equation*}
and
\begin{equation*}
  \theta_{\rho,a}(-\eta,\mu)-\theta_{\rho,a}(-\frac{\xi}{2},0)\leq \partial_1\theta_{\rho,a}(-\frac{\xi}{2},0)(-\eta+\frac{\xi}{2})
  + \frac{\|\nabla^2\theta_{\rho,a}\|_{L^\infty}}{2}
  ((\eta-\frac{\xi}{2})^2 + \mu^2),
\end{equation*}
and
\begin{equation*}
  \omega(2\eta)-\omega(\xi)\leq \omega'(\xi)(2\eta-\xi).
\end{equation*}
Then the numerator $f_{\rho,a}(\eta,\mu)\leq 2\|\nabla^2\theta_{\rho,a}\|_{L^\infty}|x_0-(\eta,\mu)|^2$, and from the polar coordinate expression we see that
the integral absolutely converges for every $\beta\in ]0,2[$.
\end{proof}

\begin{proof}[Proof of Lemma \ref{lem nLinImEs}]
In the sequel we always omit the time variable $T_*$ if there is no ambiguity.
First by the explicit formula of $u$ from Lemma 3.1 in \cite{MiaoXue}, we know
\begin{equation*}
  u(x)-u(y)\leq c(\alpha)\Big(\textrm{p.v.}\int_{\mathbb{R}^2}\frac{(x-z)^{\bot}}{|x-z|^{3+\alpha}} \theta(z)\mathrm{d}z
  -\textrm{p.v.}\int_{\mathbb{R}^2}\frac{(y-z)^{\bot}}{|y-z|^{3+\alpha}} \theta(z)\mathrm{d}z \Big),
\end{equation*}
where $c(\alpha)$ is a fixed constant. We shall split the integral into several parts. For the difference
\begin{equation*}
  \Big|\int_{|x-z|\geq 2\xi}\frac{(x-z)^{\bot}}{|x-z|^{3+\alpha}}\theta(z)\mathrm{d}z
  -\int_{|y-z|\geq 2\xi}\frac{(y-z)^\bot}{|y-z|^{3+\alpha}}\theta(z)\mathrm{d}z\Big|,
\end{equation*}
parallelling to the proof of the corresponding part in Lemma 3.2 \cite{MiaoXue} or Lemma 5.2 \cite{Kis}, we obtain that it is bounded from above by
$$
C\xi\int_{\xi}^\infty \frac{\omega(\eta)}{\eta^{2+\alpha}}\mathrm{d}\eta +C\xi^{-\alpha}\omega(\xi),
$$
with $C$ a positive constant that may depend on $\alpha$.
Then, recalling that $B_{r_0\xi}(x)$ or $B_{r_0\xi}(y)$ is the disk centered at $x$ or $y$ respectively with radius $r_0\xi$, where $r_0$ is the number introduced
in Lemma \ref{lem DispImEs}, we get
\begin{equation*}
\begin{split}
  \Big|\int_{|x-z|\leq 2\xi, z\notin B_{r_0\xi}(x)}\frac{(x-z)^\bot}{|x-z|^{3+\alpha}}\theta(z)\mathrm{d}z \Big|
  &=\Big|\int_{|x-z|\leq 2\xi, z\notin B_{r_0\xi}(x)}\frac{(x-z)^\bot}{|x-z|^{3+\alpha}}(\theta(z)-\theta(x))\mathrm{d}z \Big| \\
  & \leq C \int_{r_0\xi}^{2\xi}\frac{\omega(r)}{r^{1+\alpha}}\mathrm{d}r\leq C \xi^{-\alpha}\omega(\xi).
\end{split}
\end{equation*}
A similar estimate is true for the corresponding integral with replacing $x$ by $y$.

Now we consider the contribution of the dangerous part------the integral over $B_{r_0\xi}(x)$ and $B_{r_0\xi}(y)$. Here note that what we really need to treat is
$|(u(x)-u(y))\cdot \ell|$.
Thus from $x=\rho x_0-a$ and $y=\rho y_0-a$, we have
\begin{align}
  & \Big|\int_{B_{r_0\xi}(x)}\frac{(x-z)^\bot\cdot \ell}{|x-z|^{3+\alpha}} \theta(z)\mathrm{d}z
  -\int_{B_{r_0\xi}(y)}\frac{(y-z)^\bot\cdot\ell}{|y-z|^{3+\alpha}}\theta(z)\mathrm{d}z \Big| \nonumber\\
  = & \Big|\int_{B_{r_0\xi}(x_0)}\frac{(x_0-z)^\bot\cdot e_1}{|x_0-z|^{3+\alpha}} \theta_{\rho,a}(z)\mathrm{d}z
  -\int_{B_{r_0\xi}(y_0)}\frac{(y_0-z)^\bot\cdot e_1}{|y_0-z|^{3+\alpha}}\theta_{\rho,a}(z)\mathrm{d}z \Big|  \nonumber\\
  = & \Big| \iint_{B_{r_0\xi}(x_0)}\frac{\mu}{|x_0-(\eta,\mu)|^{3+\alpha}}\theta_{\rho,a}(\eta,\mu) \mathrm{d}\eta\mathrm{d}\mu -
  \iint_{B_{r_0\xi}(y_0)}\frac{\mu}{|y_0-(\eta,\mu)|^{3+\alpha}}\theta_{\rho,a}(\eta,\mu) \mathrm{d}\eta\mathrm{d}\mu \Big| \nonumber\\
  =& \Big|\iint_{B_{r_0\xi}(x_0)}\frac{\mu}{|x_0-(\eta,\mu)|^{3+\alpha}}\big(\theta_{\rho,a}(\eta,\mu)
  -\theta_{\rho,a}(-\eta,\mu)\big)\mathrm{d}\eta\mathrm{d}\mu\Big| \nonumber\\
  = &\Big| \iint_{B^+_{r_0\xi}(x_0)}\frac{\mu}{|x_0-(\eta,\mu)|^{3+\alpha}}\big(\theta_{\rho,a}(\eta,\mu)-\theta_{\rho,a}(-\eta,\mu)-
  \theta_{\rho,a}(\eta,-\mu)+\theta_{\rho,a}(-\eta,-\mu)\big) \mathrm{d}\eta\mathrm{d}\mu\Big|.
\end{align}
where $B_{r_0\xi}^+(x_0):= B_{r_0\xi}(x_0)\cap \{\mu>0\}$. We claim that the last expression is bounded from above by
$$-C\xi^{\beta-\alpha} \Upsilon_\beta^\bot(\xi).$$
Indeed, comparing it with \eqref{eq UpsBot}, we first
have that for every $(\eta,\mu)\in B^+_{r_0\xi}(x_0)$ and $\beta\in ]\alpha, \alpha+1]$,
\begin{equation*}
  0\leq \frac{\mu}{|x_0-(\eta,\mu)|^{3+\alpha}}= \frac{\mu}{|x_0-(\eta,\mu)|^{1+\alpha-\beta}}\frac{1}{|x_0-(\eta,\mu)|^{2+\beta}}\leq \frac{C\xi^{\beta-\alpha}}{|x_0-(\eta,\mu)|^{2+\beta}},
\end{equation*}
where $C$ is a constant that may depend on $\alpha,\beta$. Second, we get
\begin{equation*}
  |\theta_{\rho,a}(\eta,\mu)-\theta_{\rho,a}(-\eta,\mu)-
  \theta_{\rho,a}(\eta,-\mu)+\theta_{\rho,a}(-\eta,-\mu)|\leq f_{\rho,a}(\eta,\mu),
\end{equation*}
where $f_{\rho,a}(\eta,\mu)$ is defined by \eqref{eq frhoa}. In fact, it directly reduces to two obvious estimates $\theta_{\rho,a}(\eta,\mu)-\theta_{\rho,a}(-\eta,\mu)\leq \omega(2\eta)$ and $\theta_{\rho,a}(\eta,-\mu)-\theta_{\rho,a}(-\eta,-\mu)\leq \omega(2\eta)$.
Therefore, gathering the above estimates, we obtain \eqref{eq nLinImEs1}.

The proof of \eqref{eq OmeIm1} is similar to the corresponding part in Lemma 4.3 or Lemma 5.3 in \cite{Kis}, and we here sketch it.
First we consider $\omega(\xi)$ defined by \eqref{sMOC}. For the case $\xi\in ]0,\delta[$, from $\omega(\eta)\leq \frac{H}{\delta^\gamma}\eta^\gamma$ for all $\eta$, we have
\begin{equation*}
  \xi\int_\xi^\infty \frac{\omega(\eta)}{\eta^{2+\alpha}}\mathrm{d}\eta \leq \xi\int_\xi^\infty \frac{H}{\delta^{\gamma}\eta^{2+\alpha-\gamma}}\mathrm{d}\eta
  \leq \frac{1}{1+\alpha-\gamma}\frac{H}{\delta^\gamma}\xi^{\gamma-\alpha}=\frac{1}{1+\alpha-\gamma}\xi^{-\alpha}\omega(\xi).
\end{equation*}
For the case $\xi>\delta$, from $\omega(\xi)=H$, we get
\begin{equation*}
  \xi\int_\xi^\infty \frac{\omega(\eta)}{\eta^{2+\alpha}}\mathrm{d}\eta= \frac{1}{1+\alpha} H \xi^{-\alpha}= \frac{1}{1+\alpha} \xi^{-\alpha} \omega(\xi).
\end{equation*}
Then we consider $\omega(\xi,\xi_0(t))$ defined by \eqref{evtMOC}. For the case $\xi\in ]0,\xi_0(t)]$, we split the integral into two regions: from $\omega(\xi_0,\xi_0) \leq \frac{1}{1-\gamma}\omega(\xi,\xi_0)$, we know
\begin{equation*}
  \xi\int_\xi^{\xi_0} \frac{\omega(\eta,\xi_0)}{\eta^{2+\alpha}}\mathrm{d}\eta\leq \frac{1}{1+\alpha}\omega(\xi_0,\xi_0)\xi^{-\alpha}\leq \frac{1}{(1+\alpha)(1-\gamma)}\xi^{-\alpha} \omega(\xi,\xi_0);
\end{equation*}
and
\begin{equation*}
\begin{split}
  \xi\int_{\xi_0}^\infty \frac{\omega(\eta,\xi_0)}{\eta^{2+\alpha}}\mathrm{d}\eta & \leq \xi\int_{\xi_0}^\infty \frac{H}{\delta^\gamma \eta^{2+\alpha-\gamma}}\mathrm{d}\tau \leq \frac{1}{1+\alpha-\gamma}\xi \frac{H}{\delta^\gamma}\xi_0^{\gamma-\alpha-1} \\
  & \leq \frac{1}{1+\alpha-\gamma}\xi^{-\alpha}\omega(\xi_0,\xi_0)\leq \frac{1}{(1+\alpha-\gamma)(1-\gamma)}\xi^{-\alpha}\omega(\xi,\xi_0).
\end{split}
\end{equation*}
The remaining case $\xi>\xi_0(t)$ is identical to the case $\xi>0$ of $\omega(\xi)$, thus we omit it.
\end{proof}

\section{Appendix}
\setcounter{section}{6}\setcounter{equation}{0}

\subsection{Local well-posedness for $\alpha\in ]0,1[$ and $\beta\in ]2\alpha,2]$}

Our local result of \eqref{gQG} is as follows.
\begin{proposition}\label{prop local}
Let $\nu>0$, $\alpha\in]0,1[$, $\beta\in ]2\alpha,2]$ and the initial data $\theta_0\in
H^{m}$, $m>2$. Then there exists a positive constant $T$ depending only on
$\alpha,\beta,\nu,\norm{\theta_{0}}_{H^{m}}$ such that the generalized
quasi-geostrophic equation \eqref{gQG} generates a unique solution
$\theta\in \mathcal{C}([0,T],H^{m})\cap L^2([0,T],H^{m+\frac{\beta}{2}})$. Moreover we have $t^{\gamma}\theta\in
L^{\infty}(]0,T],H^{m+\gamma\beta})$ for all $\gamma\geq 0$, which implies $\theta\in \mathcal{C}^{\infty}(]0,T]\times\mathbb{R}^2)$.

Besides, we have the following blowup criterion: let $T^{*}$ be the maximal existence time of $\theta$ in
$\mathcal{C}([0,T^{*}[,H^{m})\cap L^2([0,T^*[, H^{m+\frac{\beta}{2}})$ and if $T^{*}<\infty$, then we necessarily have
\begin{equation}\label{eq blowup}
 \int_{0}^{T^{*}}\norm{\nabla \theta(t,\cdot)}_{L^{\infty}}^{2+2\alpha-\beta}\textrm{d}
 t=\infty.
\end{equation}
\end{proposition}

The proof mainly relies on the following Lemma (cf. \cite{MiaoXue})
\begin{lemma}\label{lem EstRe}
Let $v$ be a divergence free vector field over $\mathbb{R}^{n}$. For every
$q\in\mathbb{N}$, denote
\begin{equation*}
F_{q}(v,f):=S_{q+1}v\cdot\nabla\Delta_{q}f-\Delta_{q}(v\cdot\nabla
f).
\end{equation*}
Then for every $\gamma\in]0,1[$ and $p\in [1,\infty]$, there exists a positive constant $C$ such that
\begin{equation}\label{eq commutator1}
\begin{split}
 & 2^{-q\gamma}\norm{F_{q}(v,f)}_{L^p}
 \\ \leq& C\norm{\Lambda^{1-\gamma}
 v}_{L^{\infty}} \Big(\sum_{q'\leq q+4}2^{q'-q}\norm{\Delta_{q'}f}_{L^{2}} + \sum_{q'\geq q-4} 2^{(q-q')(1-\gamma)}\norm{\Delta_{q'}f}_{L^p}\Big),
\end{split}
\end{equation}
\\
Especially, in the case $n=2$ and $v=\Lambda^{\alpha}\mathcal{R}^{\bot}f$ ($\alpha\in]0,1[$), we
further have that for every $\gamma\in\big]\max\{0,\alpha\},1\big[$ and $q\in\mathbb{N}$
\begin{equation}\label{eq commutator2}
\begin{split}
 & 2^{-q\gamma}\norm{F_{q}(v,f)}_{L^{p}} \\ \leq & C \Big(\norm{\Lambda^{1-\gamma}
 v}_{L^{\infty}}\sum_{q'\geq q-4}2^{(q-q')(1-\beta)}\norm{\Delta_{q'}f}_{L^{p}}+
 \norm{\Lambda^{\alpha+1-\gamma}f}_{L^{\infty}}\sum_{|q'-q|\leq
 4}\norm{\Delta_{q'}f}_{L^{p}}\Big).
\end{split}
\end{equation}
Moreover, when $\gamma=0$, $\alpha=0$, \eqref{eq commutator1} and \eqref{eq commutator2} hold if $\norm{\Lambda^{1-\gamma}v}_{L^\infty}$ is replaed by $\norm{\nabla v}_{L^\infty}$;
and when $\gamma=1$, $\alpha=1$, then \eqref{eq commutator1} and \eqref{eq commutator2} hold if such a modification is made
$$\norm{\Lambda^{1-\gamma}v}_{L^\infty}\rightarrow\norm{v}_{B^0_{\infty,1}},\quad \norm{\Lambda^{\alpha+1-\gamma}f}_{L^\infty}\rightarrow\norm{\nabla f}_{L^\infty}.$$
\end{lemma}

Proposition \ref{prop local} is similar to Proposition 4.1-4.2 in \cite{MiaoXue}, and here we sketch the proof.
\begin{proof}[Proof of Proposition \ref{prop local}]
  Step 1: a priori estimates.

We first a priori assume that $\theta$ (and $u$) is smooth to obtain the uniform $B^m_{p,2}$ ($m>1+\frac{2}{p}$, $p\in [2,\infty[$)
estimates of $\theta$ (Note that here in the proof only $p=2$ case is used). We claim that the smooth solution $\theta(t,x)$ satisfies
\begin{equation}\label{eq BmpEst1}
  \frac{d}{dt}\norm{\theta(t)}_{B^{m}_{p,2}}^{2}+\norm{\theta(t)}_{B^{m+\frac{\beta}{2}}_{p,2}}^{2}\lesssim_{\alpha,\beta,\nu}
   \norm{ \nabla \theta}_{L^\infty}^{2+2\alpha-\beta} \norm{ \theta}_{L^\infty}^{\beta-2\alpha}\norm{\theta}_{B^{m}_{p,2}}^2
 + \norm{ \theta}_{L^{p}}^2(1+\norm{\theta}_{B^{m}_{p,2}}).
\end{equation}
Indeed, for every $q\in \mathbb{N}$, by applying the dyadic operator $\Delta_q$ to equation \eqref{gQG} we get
\begin{equation*}
 \partial_{t}\Delta_{q}\theta + S_{q+1}u
 \cdot\nabla \Delta_{q}\theta + \nu \Lambda^{\beta}\Delta_{q}\theta =
  F_{q}( u,  \theta),
\end{equation*}
where
\begin{equation*}
 F_{q}( u, \theta)=S_{q+1} u\cdot\nabla \Delta_{q} \theta- \Delta_{q}( u\cdot\nabla \theta).
\end{equation*}
Due to that $\Delta_q \theta$ is real-valued, thus multiplying both sides by $|\Delta_q \theta|^{p-2}\Delta_q\theta$ and integrating in the space variable, we obtain
\begin{equation*}
\begin{split}
 \frac{1}{p}\frac{d}{dt}\norm{\Delta_{q}
 \theta}_{L^{p}}^{p}+\nu\int_{\mathbb{R}^2} \big(\Lambda^{\beta}\Delta_q\theta(x)\big)|\Delta_{q}\theta|^{p-2}\Delta_q\theta(x)\mathrm{d}x & \leq
 \Big|\int_{\mathbb{R}^2} \big(F_q(  u, \theta)\big)(x)
 |\Delta_q\theta|^{p-2}\Delta_q\theta(x)\mathrm{d}x \Big| \\ & \leq
 \norm{F_q( u, \theta)}_{L^p} \norm{\Delta_{q}\theta}_{L^p}^{p-1}.
\end{split}
\end{equation*}
The generalized Bernstein inequality in \cite{ChenMZ} yields that an absolute constant $c>0$ (e.g. $c=1$ when $p=2$) independent of $q$ can be found such that
\begin{equation*}
  \frac{1}{p}\frac{d}{dt}\norm{\Delta_q\theta}_{L^p}^p +c\nu 2^{q\beta}\norm{\Delta_q\theta}_{L^p}^p \leq
  \norm{F_q( u, \theta)}_{L^p} \norm{\Delta_{q}\theta}_{L^p}^{p-1}.
\end{equation*}
Thus we further have
\begin{equation}\label{eq HmEs2}
  \frac{1}{2}\frac{d}{dt}\norm{\Delta_{q}
 \theta}_{L^{p}}^{2}+ \frac{c\nu}{2}2^{q\beta}
 \norm{\Delta_{q}\theta}_{L^p}^2\leq \frac{C_0}{\nu}
 \Big(2^{-q\frac{\beta}{2}}\norm{F_q( u,\theta)}_{L^p}\Big)^2.
\end{equation}
From \eqref{eq commutator2}, we know that
\begin{equation}\label{eq FqEst}
\begin{split}
 & 2^{-q\frac{\beta}{2}}\norm{F_{q}( u, \theta)}_{L^{p}}
 \\ \lesssim & \|\Lambda^{1-\frac{\beta}{2}}
  u\|_{L^{\infty}}\sum_{q'\geq
 q-4}2^{(q-q')(1-\frac{\beta}{2})}\norm{\Delta_{q'} \theta}_{L^p}+
 \|\Lambda^{\alpha+1-\frac{\beta}{2}} \theta\|_{L^{\infty}}\sum_{|q'-q|\leq
 4}\norm{\Delta_{q'} \theta}_{L^p}
\end{split}
\end{equation}
Also notice that for some number $K\in\mathbb{N}$
\begin{equation*}
 \begin{split}
   \|\Lambda^{1-\frac{\beta}{2}} u\|_{L^{\infty}}+ \|\Lambda^{\alpha+1-\frac{\beta}{2}} \theta\|_{L^{\infty}} &\lesssim
   \|\Lambda^{1-\frac{\beta}{2}} \Lambda^{\alpha}\mathcal{R}^{\bot}\theta\|_{\dot{B}^{0}_{\infty,1}}+
   \|\Lambda^{\alpha+1-\frac{\beta}{2}} \theta\|_{\dot{B}^{0}_{\infty,1}}   \\  & \lesssim
   \sum_{k=-\infty}^{K -1}
   2^{k(\alpha+1-\frac{\beta}{2})}\|\dot\Delta_k \theta\|_{L^{\infty}}
   + \sum_{k= K }^{\infty}2^{-k(\frac{\beta}{2}-\alpha)}\|\dot\Delta_k \nabla\theta\|_{L^{\infty}}   \\
   &\lesssim 2^{K(\alpha+1-\frac{\beta}{2})} \| \theta\|_{L^{\infty}}+ 2^{K(\alpha-\frac{\beta}{2})}
   \norm{ \nabla\theta}_{L^\infty},
 \end{split}
\end{equation*}
thus choosing $K$ satisfying $\norm{\theta}_{L^\infty}2^{K}\thickapprox \norm{ \nabla \theta}_{L^\infty}$, we infer
\begin{equation}\label{eq simplInq}
 \|\Lambda^{1-\frac{\beta}{2}} u\|_{L^{\infty}}+ \|\Lambda^{\alpha+1-\frac{\beta}{2}} \theta\|_{L^{\infty}}
 \lesssim \norm{  \nabla \theta}_{L^\infty}^{1+\alpha-\frac{\beta}{2}}\norm{ \theta}_{L^\infty}^{\frac{\beta}{2}-\alpha}.
\end{equation}
Plunging the above two estimates \eqref{eq simplInq} and \eqref{eq FqEst} into inequality \eqref{eq HmEs2},
then multiplying both sides by $2^{2qm}$ and summing up over $q\in\mathbb{N}$, from Young inequality we obtain
\begin{equation}\label{eq hmEs-HF}
 \frac{1}{2}\frac{d}{dt}\sum_{q\in\mathbb{N}}2^{2qm}\norm{\Delta_{q}
 \theta}_{L^p}^2+\frac{c\nu}{2}\sum_{q\in\mathbb{N}}2^{2q(m+\frac{\beta}{2})}
 \norm{\Delta_q\theta}_{L^p}^{2}
 \lesssim \frac{1}{\nu} \norm{  \nabla \theta}_{L^\infty}^{2+2\alpha-\beta} \norm{
 \theta}_{L^\infty}^{\beta-2\alpha} \norm{\theta}_{B^m_{p,2}}^2.
\end{equation}
On the other hand, we apply the low frequency operator $\Delta_{-1}$ to \eqref{gQG} to get
\begin{equation*}
 \partial_t\Delta_{-1}\theta  =- \nu \Lambda^\beta\Delta_{-1}\theta-\Delta_{-1}\big(  u \cdot\nabla  \theta \big).
\end{equation*}
Multiplying both sides by $|\Delta_{-1}\theta|^{p-2}\Delta_{-1}\theta$, integrating over the spatial variable and using the positivity formula of the dissipative term, we obtain
\begin{equation*}
\begin{split}
 \frac{1}{p}\frac{d}{dt}\norm{\Delta_{-1}\theta}_{L^p}^p
  & \leq \Big|\int_{\mathbb{R}^2} \mathrm{div}\Delta_{-1}\big(   u \, \theta\big)(x) \,
 |\Delta_{-1}\theta|^{p-2}\Delta_{-1}\theta(x)\mathrm{d}x\Big| \\
 & \lesssim \norm{  u}_{L^\infty}\norm{  \theta}_{L^p}\|\Delta_{-1}\theta\|_{L^p}^{p-1}.
\end{split}
\end{equation*}
We see that
\begin{equation}\label{eq fact1}
\begin{split}
 \norm{  u}_{L^\infty} & \leq \Big(\sum_{j\leq -1}+ \sum_{j\geq 0}\Big) \|\dot \Delta_j
 \Lambda^{\alpha}\mathcal{R}^{\bot} \theta\|_{L^\infty} \\
 & \lesssim\sum_{j\leq -1}2^{j\alpha}\|\dot\Delta_j   \theta\|_{L^\infty}+ \sum_{j\geq 0} 2^{j(\alpha-1)} \|\dot \Delta_j
 \nabla \theta\|_{L^\infty} \lesssim \norm{ \theta}_{L^\infty}+  \norm{\nabla\theta}_{L^\infty},
\end{split}
\end{equation}
thus we have
\begin{equation}\label{eq hmEs-LF}
 \frac{1}{2}\frac{d}{dt}\norm{\Delta_{-1}\theta}_{L^p}^2+ \frac{c\nu}{2} 2^{-\beta}\norm{\Delta_{-1} \theta}_{L^p}^2
 \lesssim (\norm{ \theta}_{B^m_{p,2}}+\nu)\norm{ \theta}_{L^p}^2 .
\end{equation}
Multiplying \eqref{eq hmEs-LF} by $2^{-2m}$ and combining it with \eqref{eq hmEs-HF} leads to \eqref{eq BmpEst1}.

Set $G(t):= \norm{\theta(t)}_{B^m_{p,2}}^2+\int_0^t \norm{\theta(\tau)}^2_{B^{m+\beta/2}_{p,2}}\mathrm{d}\tau$. Then \eqref{eq BmpEst1} reduces to
\begin{equation*}
  \frac{d}{dt} G(t)\leq C (G(t)^2+  G(t)),
\end{equation*}
where $C$ will depend on $\nu,\alpha,\beta,m$. Gronwall inequality ensures that for every
\begin{equation}\label{eq timBd}
T<\frac{1}{C}\log (1+ 1/\|\theta_0\|_{B^m_{p,2}}^2),
\end{equation}
we have
\begin{equation}\label{eq BmpUniBd}
  \sup_{t\in [0,T]}\norm{\theta(t)}_{B^m_{p,2}}^2 +\int_0^T \norm{\theta(\tau)}^2_{B^{m+\beta/2}_{p,2}}\mathrm{d}\tau \leq
  \frac{\norm{\theta_0}_{B^m_{p,2}}^2}{(\|\theta_0\|_{B^m_{p,2}}^2+1)e^{-CT}-\|\theta_0\|_{B^m_{p,2}}^2}.
\end{equation}
Note that \eqref{eq timBd} and \eqref{eq BmpUniBd} hold for all $G(0)\leq \|\theta_0\|_{B^m_{p,2}}^2$.

Step 2: Uniqueness

Let $\theta^{(1)}$, $\theta^{(2)}\in L^\infty ([0,T], H^{m}(\mathbb{R}^2))$, $m>2$ be two smooth solutions to the
generalized quasi-geostrophic equation \eqref{gQG} with the same
initial data. Denote $u^{(i)}=\Lambda^{\alpha}\mathcal{R}^{\bot}\theta^{(i)}$,
$i=1,2$, $\delta \theta=\theta^{(1)}-\theta^{(2)} $, $\delta u=
u^{(1)}-u^{(1)}$, then we write the difference equation as
\begin{equation*}
  \partial_t\delta\theta + u^{(1)}\cdot\nabla\delta\theta + \nu
  \Lambda^\beta\delta\theta = -\delta u\cdot\nabla \theta^{(2)},\quad
  \delta\theta|_{t=0}=0
\end{equation*}
By $L^2$ energy method, we have
\begin{equation*}
\begin{split}
  \frac{1}{2}\frac{d}{dt}\norm{\delta\theta(t)}_{L^2}^2 + \nu \|\Lambda^{\frac{\beta}{2}}\delta\theta(t)\|_{L^2}^2 & \leq
  \Big|\int_{\mathbb{R}^2} (\delta u\cdot \nabla\theta^{(2)})(t,x) \;\delta\theta(t,x)\mathrm{d}x\Big| \\
  & \leq \|\delta u\cdot\nabla \theta^{(2)}(t)\|_{\dot H^{-\frac{\beta}{2}}} \|\delta\theta(t)\|_{\dot H^{\frac{\beta}{2}}} \\
  & \leq C_{\nu,\alpha,\beta}\|\Lambda^\alpha\mathcal{R}^\bot \delta\theta(t)\|^2_{\dot H^{-\alpha}}\|\nabla\theta^{(2)}(t)\|_{\dot H^{\alpha+1-\beta/2}}^2
  +\frac{\nu}{2} \|\Lambda^{\frac{\beta}{2}}\delta\theta(t)\|_{L^2}^2,
\end{split}
\end{equation*}
where in the last line we have used Young inequality and the following classical estimate that for every divergence-free $f$ and every $s,t<1$, $s+t>-1$,
\begin{equation*}
 \norm{f\cdot\nabla g}_{\dot H^{s+t-1}}\lesssim_{s,t} \norm{f}_{\dot H^s}\norm{\nabla g}_{\dot H^t}.
\end{equation*}
From $H^m\hookrightarrow \dot H^{2+\alpha-\beta/2}$ continuously, we further obtain
\begin{equation*}
  \frac{d}{dt}\norm{\delta\theta(t)}_{L^2}^2 \leq C_{\nu,\alpha,\beta} \norm{\delta\theta(t)}_{L^2}^2 \|\theta^{(2)}(t)\|_{H^m}^2,
\end{equation*}
thus the Gronwall inequality leads to $\delta\theta(t)\equiv 0$ for all $t\in [0,T]$.

Step 3: Existence

We first regulate the system \eqref{gQG} to get
\begin{equation}\label{Re-gQG}
 \begin{cases}
  \begin{split}
    &\theta^{N}_{t}+J_{N}
   \big( u^N\cdot\nabla\theta^N\big)
   +\nu \Lambda^{\beta}\theta^{N} =0 \\
    &u^N= \Lambda^\alpha\mathcal{R}^{\perp}\theta^N, \quad
    \theta^N|_{t=0}=J_N\theta_{0},
  \end{split}
 \end{cases}
\end{equation}
where $J_N:L^2\rightarrow J_N L^2$, $N\in\mathbb{Z}^+$ is the Friedrich projection operator such that $\widehat{J_N f}(\zeta)= 1_{B_N}(\zeta)\hat{f}(\zeta)$.
By Cauchy-Lipschitz theorem, for every $N\in\mathbb{Z}^+$ there exists a unique global solution $\theta^N\in \mathcal{C}^1([0,\infty[,H^\infty(\mathbb{R}^2))$ to the regularized system \eqref{Re-gQG}. Then almost paralleling to the proof in the step 1, we know the uniform estimate that for all $T<\frac{1}{C}\log (1+ 1/\|\theta_0\|_{H^m}^2)$
\begin{equation}\label{eq HmUniBd}
  \sup_{t\in [0,T]}\|\theta^N(t)\|_{H^m}^2 +\int_0^T \|\theta^N(\tau)\|^2_{H^{m+\beta/2}}\mathrm{d}\tau \leq
  \frac{\norm{\theta_0}_{H^m}^2}{(\|\theta_0\|_{H^m}^2+1)e^{-CT}-\|\theta_0\|_{H^m}^2},
\end{equation}
where $C$ is a positive constant depending only on $\nu,\alpha,\beta,m$. From the uniform estimate \eqref{eq HmUniBd} and the uniqueness result, similarly as treating the corresponding part in \cite{MiaoXue}, we can prove that $\theta^N$ is a Cauchy sequence in $\mathcal{C}([0,T];L^2(\mathbb{R}^2))$, which implies a strong convergence to a function $\theta\in\mathcal{C}([0,T];L^2)$. By a classical method we know that $\theta$ is a solution of the limiting system \eqref{gQG},
and satisfies $\theta\in L^\infty([0,T]; H^m(\mathbb{R}^2))\cap L^2([0,T];H^{m+\beta/2}(\mathbb{R}^2))$. Moreover, similarly as \cite{MiaoXue}, we can prove the issues of the time continuity in $H^m$, the $\mathcal{C}^\infty$ smoothness in $]0,T]\times \mathbb{R}^2$ and the blowup criterion of the solution.

\end{proof}

\subsection{Global existence of weak solutions for $\alpha\in ]0,1[$ and $\beta\in ]2\alpha,2]$}

The main result in this subsection is as follows.
\begin{proposition}\label{prop GlWeak}
Let $\nu>0$, $\alpha\in ]0,1[$, $\beta\in ]2\alpha,2]$ and $\theta_0\in L^2(\mathbb{R}^2)$. Then there exists a global weak solution $\theta \in L^\infty([0,\infty[; L^2(\mathbb{R}^2))\cap L^2([0,\infty[;\dot H^{\beta/2}(\mathbb{R}^2))$ for the generalized quasi-geostrophic equation \eqref{gQG}. Moreover, $\theta$ will satisfy the following energy inequality
\begin{equation}\label{eq EnegEs}
  \|\theta(t)\|_{L^2(\mathbb{R}^2)}^2 + 2\nu \int_0^t\|\theta (\tau)\|^2_{\dot H^{\frac{\beta}{2}}(\mathbb{R}^2)}\mathrm{d}\tau \leq \|\theta_0\|_{L^2(\mathbb{R}^2)}^2,\quad t>0.
\end{equation}
\end{proposition}

The proof follows from the standard process of establishing weak solutions (cf. \cite{Tem,Lemarie,FriV}), and here we sketch it.
\begin{proof}[Proof of Proposition \ref{prop GlWeak}]
We consider the following approximate system
\begin{equation}\label{eq APgQG}
  \begin{cases}
    \partial_t\theta^\epsilon + u^\epsilon \cdot\nabla \theta^\epsilon +\nu\Lambda^\beta \theta^\epsilon -\epsilon\Delta \theta^\epsilon=0 \\
    u^\epsilon= \Lambda^\alpha \mathcal{R}^{\bot}\theta^\epsilon, \quad \theta^\epsilon|_{t=0}= \psi_\epsilon*\theta_0,
  \end{cases}
\end{equation}
where $\psi(x)\in \mathcal{C}^\infty_c(\mathbb{R}^2)$ is a radial positive function satisfying $\int \psi =1$ and $\psi_\epsilon(x)=\epsilon^{-2}\psi(x/\epsilon)$. Clearly, $\|\psi_\epsilon*\theta_0\|_{H^m}\lesssim_\epsilon \|\theta_0\|_{L^2}$ for all $m>0$. Thus for $m>2$ and $\epsilon>0$, from Theorem \ref{thm1}, the approximate generalized quasi-geostrophic equation \eqref{eq APgQG} exists a unique global smooth solution
$$
\theta^\epsilon\in \mathcal{C}([0,\infty[;H^{m}) \cap \mathcal{C}^{\infty}(]0,\infty[\times \mathbb{R}^{2}).
$$
Besides, from $\|\psi_\epsilon*\theta_0\|_{L^2}\leq \|\theta_0\|_{L^2}$, we also have the uniform energy estimate
\begin{equation}\label{eq EnegEs1}
   \|\theta^\epsilon(T)\|_{L^2(\mathbb{R}^2)}^2 + 2\nu \int_0^T\|\theta^\epsilon (\tau)\|^2_{\dot H^{\frac{\beta}{2}}(\mathbb{R}^2)}\mathrm{d}\tau \leq \|\theta_0\|_{L^2(\mathbb{R}^2)}^2,\quad \forall T>0,
\end{equation}
that is, $\theta^\epsilon$ is uniformly bounded in $\mathcal{C}([0,T]; L^2(\mathbb{R}^2))\cap L^2([0,T]; \dot H^{\beta/2}(\mathbb{R}^2))$ for every $T>0$. This ensures that, up to the subsequence, $\theta^\epsilon$ converges weakly to some function $\theta\in L^\infty([0,T]; L^2)\cap L^2([0,T]; \dot H^{\beta/2})$ (in $L^\infty([0,T];L^2)$ the convergence is weak-$*$). But it is not sufficient to pass to the limit of the nonlinear term of \eqref{eq APgQG} in the weak framework.

We further claim that,
\begin{equation}\label{eq CovSL2}
  \theta^\epsilon\rightarrow \theta\quad \textrm{strongly in}\; L^2([0,T]; L_{\mathrm{loc}}^2(\mathbb{R}^2)).
\end{equation}
We shall use the classical Aubin-Lions compactness Lemma (cf. \cite{Tem}) to prove it. For any compact subset $\mathcal{O}\subset \mathbb{R}^2$,
since the mapping $f \rightarrow \chi_{\mathcal{O}}f$ is
compact from $H^{\beta/2}$ to $L^2$ (cf. \cite{Lemarie}), with $\chi_{\mathcal{O}}(x)\in \mathcal{C}^\infty_0 (\mathbb{R}^2)$ satisfying that
it is supported in a compact subset $\mathcal{O}'$ and $\chi_{\mathcal{O}}(x)=1$ for all $x\in\mathcal{O}$,
we know that the sequence $\{ \chi_{\mathcal{O}}\theta^\epsilon\}_{\epsilon>0}$ is compact in
$ L^2$. Thus to guarantee \eqref{eq CovSL2} it suffices to find some suitable Banach space $X$ and $a_0\in ]1,\infty[$ such that $L^2(\mathcal{O}')\hookrightarrow X(\mathcal{O}') $ continuously and $\partial_t\theta^\epsilon$ uniformly bounded in $ L^{a_0}([0,T]; X(\mathcal{O}'))$. In fact, we shall prove
\begin{equation}\label{eq CovEs2}
  \partial_t \theta^\epsilon\; \textrm{is uniformly}\; \textrm{bounded in}\; L^{\frac{4}{3}}([0,T]; W^{-2, \frac{8}{8-\beta}}(\mathcal{O}')).
\end{equation}
From \eqref{eq EnegEs1}, interpolation and Sobolev embedding, we have $\theta^{\epsilon}\in L^4 ([0,T];L^{\frac{8}{4-\beta}}(\mathbb{R}^2)$.
On the other hand, due to that the Riesz transform is bounded in $L^2(\mathbb{R}^2)$ and $H^{\beta/2}(\mathbb{R}^2)\hookrightarrow \dot H^{\alpha}(\mathbb{R}^2)$ continuously, we immediately get $u^\epsilon\in L^2([0,T]; L^2(\mathbb{R}^2))$ uniformly in $\epsilon$. Thus $\mathrm{div}(u^\epsilon \theta^\epsilon)$ is uniformly bounded in $L^{\frac{4}{3}}([0,T]; W^{-1, \frac{8}{8-\beta}}(\mathbb{R}^2))$. For the dissipative terms, we get
$\epsilon\Delta\theta^\epsilon$ is uniformly bounded in $L^2{[0,T]; H^{\frac{\beta}{2}-2}(\mathbb{R}^2)} $, and $\nu\Lambda^\beta\theta^\epsilon$
is uniformly bounded in $L^2([0,T]; H^{-1}(\mathbb{R}^2))$. Hence, if the spacial variable is restricted in $\mathcal{O}'$, we obtain $\partial_t \theta^\epsilon$ is uniformly bounded in $L^{\frac{4}{3}}([0,T]; W^{-1, \frac{8}{8-\beta}}(\mathcal{O}')) + L^2([0,T]; H^{\frac{\beta}{2}-2}(\mathcal{O}')) + L^2([0,T]; H^{-1}(\mathcal{O}'))$; and thus uniformly in $L^{\frac{4}{3}}([0,T]; W^{-2, \frac{8}{8-\beta}}(\mathcal{O}'))$ from the continuous embedding.

Moreover, we also have
\begin{equation}\label{eq uCovSL2}
u^\epsilon\rightarrow u\;\textrm{strongly in}\; L^2([0,T];L^2_{\mathrm{loc}}(\mathbb{R}^2)),
\end{equation}
that is,
$$
\Lambda^\alpha\mathcal{R}^{\bot} \theta^\epsilon\rightarrow \Lambda^\alpha\mathcal{R}^{\bot}\theta\;\textrm{strongly in}\; L^2([0,T];L^2_{\mathrm{loc}}(\mathbb{R}^2)).
$$
The proof is similar to \eqref{eq CovSL2}. From \eqref{eq EnegEs}, we know that for all $T>0$, $u^\epsilon$ is uniformly bounded in $L^\infty([0,T]; \dot H^{-\alpha}(\mathbb{R}^2))\cap L^2([0,T]; \dot H^{\frac{\beta}{2}-\alpha})$. By interpolation, $u^\epsilon\in L^2([0,T]; L^2(\mathbb{R}^2))$ uniformly in $\epsilon$, thus $u^\epsilon\in L^2([0,T]; H^{\beta/2-\alpha}(\mathbb{R}^2))$ uniformly in $\epsilon$.
For any compact subset $\mathcal{O}\subset\mathbb{R}^2$, since the sequence
$\{ \chi_{\mathcal{O}}u^\epsilon\}_{\epsilon>0}$ is compact in $ L^2$,
it remains to show that $\partial_t u^\epsilon$ is uniformly bounded in $ L^{\frac{4}{3}}([0,T]; W^{-2,\frac{8}{8-\beta}}(\mathcal{O}'))$.
Indeed, from the equation
\begin{equation*}
  \partial_t\Lambda^\alpha\mathcal{R}^{\bot}\theta^\epsilon = - \Lambda^\alpha\mathcal{R}^{\bot}\mathrm{div}(u^\epsilon\theta^\epsilon) -\nu
  \Lambda^{\alpha+\beta}\mathcal{R}^{\bot}\theta^\epsilon - \epsilon\Lambda^{\alpha+2}\mathcal{R}^{\bot}\theta^\epsilon,
\end{equation*}
we have
\begin{equation*}
\begin{split}
  \|\partial_t u^\epsilon\|_{L^{4/3}_T W^{-2, \frac{8}{8-\beta}}(\mathcal{O}')}& \lesssim \|\mathrm{div}(u^\epsilon\theta^\epsilon) \|_{L^{4/3}_T W^{-1, \frac{8}{8-\beta}}(\mathbb{R}^2)} + \nu \| \theta^\epsilon\|_{L^2_T H^{\alpha+\beta-2}(\mathbb{R}^2)} + \| \theta^\epsilon\|_{L^2_T H^{\alpha}(\mathbb{R}^2)}\\
  & \lesssim_\nu \|u^\epsilon\theta^\epsilon \|_{L^{4/3}_T L^{ \frac{8}{8-\beta}}(\mathbb{R}^2)} + \| \theta^\epsilon\|_{L^2_T H^{\beta/2}(\mathbb{R}^2)}.
\end{split}
\end{equation*}
Then \eqref{eq uCovSL2} follows from the classical Aubin-Lions Lemma.

Based on \eqref{eq CovSL2} and \eqref{eq uCovSL2}, we can send to the limit in \eqref{eq APgQG}. Indeed, for any $\phi\in\mathcal{C}^\infty_c (]0,\infty[\times\mathbb{R}^2)$, we have
\begin{equation*}
\begin{split}
  \Big|\int (u^\epsilon\theta^\epsilon)\cdot\nabla \phi -(u\theta)\cdot\nabla\phi\Big|  \leq &  \Big|\int \theta^\epsilon(u^\epsilon-u)\cdot \nabla \phi\Big| + \Big|\int (\theta^\epsilon -\theta)u\cdot\nabla\phi\Big| \\
   \leq &  \|u^\epsilon-u\|_{L^2_T L^2(\mathcal{O})}\|\theta^\epsilon\|_{L^\infty_T L^2}\|\phi\|_{L^2_T W^{1,\infty}} \\
   & + \|\theta^\epsilon-\theta\|_{L^2_T L^2(\mathcal{O})}\|u\|_{L^2_T L^2} \|\phi\|_{L^\infty_T W^{1,\infty}}.
\end{split}
\end{equation*}

Finally, we show that $\theta$ is the weak solution of the generalized quasi-geostrophic equation \eqref{gQG}, and \eqref{eq EnegEs} follows
from \eqref{eq EnegEs1} by a limiting argument (cf. \cite{Lemarie}).

\end{proof}




\begin{thebibliography}{60}
{\small
\bibitem{Caffarelli}L. Caffarelli and  V. Vasseur, Drift
                    diffusion equations with fractional diffusion and the
                    quasi-geostrophic equations. Arxiv,
                    math.AP/0608447, To appear in Annals of Math.
\bibitem{CafS}L. Caffarelli and L. Silvestre. An extension problem related to the fractional Laplacian.
                  Communications in PDE, \textbf{32}, Issue 8(2007), 1245-1260.
\bibitem{ChenMZ}Q. Chen, C. Miao and Z. Zhang, A new Bernstein's inequality
 and the 2D dissipative quasi-geostrophic equation. Comm. Math.
 Phys. \textbf{271}(2007), 821-838.
\bibitem{Con3}  P. Constantin, D. Cordoba  and J. Wu,
           On the critical dissipative quasi-geostrophic equation.
           Indiana Univ. Math. J.  \textbf{50}(2001), 97-107.
\bibitem{Con1} P. Constantin, A.J. Majda and E. Tabak,
               Formation of strong fronts in the 2-D quasigeostrophic thermal active scalar. Nonlinearity \textbf{7}(1994), 1495-1533.
\bibitem{Con2}  P. Constantin  and J. Wu, Behavior of solutions of 2D quasi-geostrophic equations.  SIAM J. Math. Anal. \textbf{30}(1999), 937-948.
\bibitem{ConW}P. Constantin and J. Wu, Regularity of H\"older continuous solutions of the supercritical quasi-geostrophic equation,
          Ann. Inst. H. Poincare Anal. Non Lineaire, \textbf{25}(2008), No.6, 1103-1110.
\bibitem{ConIW}P. Constantin, G. Iyer and J. Wu, Global regularity for a modified critical dissipative quasi-geostrophic equation.
                      Indiana University Mathematics Journal, \textbf{57}(2008), 2681-2692.
%
\bibitem{CorC}A. C\'ordoba and D. C\'ordoba, A maximum principle applied to
 the quasi-geostrophic equations. Comm. Math. Phys. \textbf{249}(2004),  511-528
\bibitem{Dab}M. Dabkowski, Eventual regularity of the solutions to the supercritical dissipative quasi-geostrophic equation,
             arXiv:math/1007.2970
%
\bibitem{Dong1}H. Dong, Well-posedness for a transport equation with nonlocal velocity. J. Funct. Anal., \textbf{255} (2008), no. 11, 3070-3097.
%
\bibitem{FriV}S. Friedlander and V. Vicol, Global well-posedness for an advection-diffusion equation arising in magneto-geostrophic dynamics.
                 Arxiv, math.AP/1007.1211v1.
\bibitem{HmiK}T. Hmidi and S. Keraani, Global solutions of the
 supercritical 2D dissipative quasi-geostrophic equation,
 \textit{Adv. Math.} \textbf{214}(2007), 618-638
\bibitem{KisN}A. Kiselev and F. Nazarov, A variation on a theme of Caffarelli and Vasseur, Zap. Nauchn. Sem. POMI
            \textbf{370}(2010), 58-72.
\bibitem{KisNV}A. Kiselev, F. Nazarov and A. Volberg, Global well-posedness
 for the critical 2D dissipative quasi-geostrophic equation, Invent. Math. \textbf{167}(2007), 445-453.
\bibitem{Kis}A. Kiselev, Nonlocal maximum principle for active scalars. Arxiv, math.AP/1009.0542.

\bibitem{Lemarie}P.G. Lemari\'{e}-Rieusset, Recent developments in the Navier-Stokes problem, Chapman and Hall/CRC, 2002.
\bibitem{MiaoXue}C. Miao and L. Xue, Global well-posedness for a modified critical dissipative quasi-qeostrophic equation. Arxiv, math.AP/0901.1368.
\bibitem{Resnick}S. Resnick, Dynamical problems in nonlinear advective partial differential equations,
            Ph.D. thesis, University of Chicago, 1995.
\bibitem{Sil}L. Silvestre, Eventual regularization for the slightly supercritical quasi-geostrophic equation, Ann. Inst. H.
       Poincar\'e Anal. Non Lin\'eaire \textbf{27}(2010), 693-704.
\bibitem{Tem}R. Temam, Navier-Stokes equations, North Holland, 1977.
\bibitem{Wu 01}J. Wu, Dissipative quasi-geostrophic equations with $L^p$ data, Electron J. Differential Equations,  \textbf{2001}(2001), 1-13.
\bibitem{Wu 04} J. Wu,  Global solutions of the 2D dissipative quasi-geostrophic equation in Besov spaces.
              SIAM J. Math. Anal. \textbf{36}(2004), 1014-1030.
}
\end{thebibliography}
\end{document}